\documentclass[a4paper]{article}

\usepackage{changebar}
\usepackage{amssymb,amsmath, amsfonts}
\usepackage{caption}
\usepackage{authblk}
\usepackage{graphicx}
\usepackage{url}
\usepackage{lineno}

\def\otm{\otimes}
\def\str{\rightarrow}

\def\mj{\mbox{\bf 1}}
\def\mjs{\mbox{\scriptsize\bf 1}}

\def\prop#1#2{\vspace{2ex} \noindent{\sc #1.} {\it #2} \par \vspace{2ex}}
\def\dkz{\noindent{\sc Proof. }}
\def\qed{\hfill $\dashv$\vspace{2ex}}

\def\mn{\!-\!}
\def\strt{\stackrel{\textbf{.}\,}{\rightarrow}}
\def\oW{\overline{W}}

\def\cM{{\cal M}}
\def\WM{\oW\!\cM}

\def\vmj{\vec{\mj}}
\def\le{\left(\!\!\!}
\def\re{\!\!\!\right)}

\begin{document}

\title{The $n$-fold reduced bar construction}
\author{Sonja Lj. \v{C}uki\'c}
\author{Zoran Petri\'c}
\affil{Mathematical Institute, SANU,\\ Knez Mihailova 36, p.f.\
367,\\ 11001 Belgrade, Serbia\\ \texttt{sonja@mi.sanu.ac.rs},
\texttt{zpetric@mi.sanu.ac.rs} }

\date{}
\maketitle

\vspace{-3ex}

\begin{abstract}
\noindent This paper is about a correspondence between monoidal
structures in categories and $n$-fold loop spaces. We developed a
new syntactical technique whose role is to substitute the
coherence results, which were the main ingredients in the proof
that the Segal-Thomason bar construction provides an appropriate
simplicial space. The results we present here enable more common
categories to enter this delooping machine. For example, such as
the category of finite sets with two monoidal structures brought
by the disjoint union and Cartesian product.

\end{abstract}

\vspace{.3cm}


\vspace{.5ex}

\noindent {\small {\it Keywords$\,$}: bar construction, monoidal
categories, infinite loop spaces}

\vspace{.5ex}

\noindent {\small {\it Acknowledgements$\,$}: This work was
supported by projects of the Ministry of Education, Science, and
Technological Development of the Republic of Serbia (ON174020 and
ON174026). }

\section{Introduction}

A correspondence between monoidal structures in categories and
loop spaces was initially established by Stasheff in \cite{S63}.
Since then, a connection of various algebraic structures on a
category with 1-fold, 2-fold, $n$-fold, and infinite loop spaces
is a subject of many papers (see \cite{ML63}, \cite{S74},
\cite{M74}, \cite{T79}, \cite{JS93}, \cite{D94}, \cite{McCS},
\cite{BFSV}, and references therein). The categories in question
are usually equipped with one or several monoidal structures, and
natural transformations providing symmetry, braiding, or some
other kind of interchange between these structures. There are two
main approaches to the subject. One is operadic and the other is
through the Segal-Thomason bar construction, which we simply call
\emph{reduced} bar construction, as in \cite{T79}. The latter, to
which we will keep to throughout the paper, is an approach to the
Quillen plus construction and it is the initial step connecting
various monoidal categories with loop spaces.

The $n$-fold reduced bar construction based on an $n$-fold
monoidal category $\cM$ is an iteration of a construction of a
simplicial object based on a monoid in a category whose monoidal
structure is given by finite products. The goal is to obtain a lax
functor $\WM$ from $(\Delta^{{\!op}})^n$, the $n$th power of the
opposite of the simplicial category, to the category \emph{Cat},
of categories and functors, such that
\[
\WM(k_1,\ldots,k_n)=\cM^{k_1\cdot\ldots\cdot k_n}.
\]

The main result of this paper states that for every $n\geq 2$, the
$n$-fold reduced bar construction delivers a certain lax functor.
This is what we mean by \emph{correctness} of the reduced bar
construction. We prove this result gradually---the cases $n=2$,
$n=3$ and $n\geq 3$ are dealt with respectively in Theorems~4.5,
6.5 and 8.5.

Following the ideas of \cite[Section~2]{BFSV}, we show in
Section~9, that the lax functor $\WM$ satisfies some additional
conditions. Roughly speaking, some particular arrows of
$(\Delta^{{\!op}})^n$, which are built out of face maps
corresponding to projections, have to be mapped by $\WM$ to
identities. Such a lax functor is called Segal's in \cite{P14}.

By applying Street's rectification to $\WM$ (see \cite{S72}) one
obtains a functor $V$, with the same source and target as $\WM$.
From \cite[Corollary~4.4]{P14}, when $B$ is the classifying space
functor, it follows that $B\circ V$ is a multisimplicial space
with some properties guaranteing that, up to group completion (see
\cite{S74} and \cite{McDS76}), the realization of this
multisimplicial space is an $n$-fold delooping of~$B\cM$ (see
\cite[Theorem~5.1]{P14}). A thorough survey of results concerning
these matters is given in \cite{P14} and the case $n=2$ is
considered separately in \cite{P15}.

\vspace{.3ex}

This paper is strongly influenced by \cite{BFSV}. One can find the
main ideas followed by us in Sections 0, 1 and 2 of that paper.
Also, the reader should consider \cite{F95} as an earlier source
of these ideas. The notions of two, three and $n$-fold monoidal
categories used in \cite{BFSV} and the corresponding notions used
in this paper are compared in Sections~2, 5 and 7. The case $n=2$
is studied systematically in Section~2.

A definition of $n$-fold monoidal category is usually inductive as
one starts with the 2-category \emph{Cat} whose monoidal structure
is given by 2-products. The 0-cells of a 2-category $Mon(Cat)$ are
pseudomonoids (or monoids) in \emph{Cat}, i.e.,\ monoidal (or
strict monoidal) categories. Then one makes a choice what to
consider to be the 1-cells of $Mon(Cat)$, i.e.,\ how strictly they
should preserve the monoidal structure. The monoidal structure of
$Mon(Cat)$ is again given by 2-products. A pseudomonoid (or a
monoid) in $Mon(Cat)$ is a (strict) two-fold monoidal category and
if we iterate this procedure with the same degree of strictness,
we obtain one possible notion of $n$-fold monoidal category.

Joyal and Street, \cite{JS93}, dealt with such a concept having in
its basis a certain 2-category of monoidal categories, strong
monoidal functors, and monoidal transformations. They showed that
such a degree of strictness leads to a sequence of categorial
structures starting with monoidal categories, then the braided
monoidal categories as the two-fold monoidal categories and
symmetric monoidal categories as the n-fold monoidal categories
for $n\geq 3$.

Balteanu, Fiedorowicz, Schw\" anzl, and Vogt, \cite{BFSV},
considered a variant of $Mon(Cat)$ consisting of strict monoidal
categories, monoidal functors (in which the interchange between
multiplicative structures need not be invertible), and monoidal
transformations. This was an important advance leading to a
definition of $n$-fold monoidal categories without stabilization
at $n=3$. However, they did not laxify the appropriate
interchanges for units, which were treated in their work as strict
as possible.


The reason to stop at that notion is probably the impossibility of
proving an appropriate coherence result for the completely
balanced notion of iterated monoidal categories. Monoidal units
usually produce difficulties in coherence results (cf.
\cite{KML71}). The situation brought by diversifying monoidal
units in the case of $n$-fold monoidal structures is very
complicated.


The idea of \cite{DPInt} and \cite{PT13} was to laxify the
interchanges for units as much as the coherence allows. Trimble
and the second author showed that a coherence result for
pseudocommutative pseudomonoids, for which some structural
constraints are invertible, in a 2-category of symmetric monoidal
categories, lax symmetric monoidal functors, and monoidal
transformations is sufficient for the reduced bar construction.

In this paper we consider the variant of $Mon(Cat)$ in which the
interchange between multiplicative structures and interchange
between units need not be invertible, i.e.,\ a 2-monoidal category
of monoidal categories, lax monoidal functors, and monoidal
natural transformations. This is the basis used by Aguiar and
Mahajan, \cite{AM10}, for the definition of the notion of $n$-fold
monoidal category. The possibility of defining $n$-fold monoidal
structures with respect to such a basis is much less explored,
perhaps because of difficulties in proving the corresponding
coherence results. Such a coherence result usually guarantees
commutativity of all the diagrams in $n$-fold monoidal categories
relevant for the reduced bar construction.


Our result is not of the form to prove the coherence and not to
worry about the lax conditions. We have developed a syntactical
technique whose role is to substitute the coherence results. The
correctness of the reduced bar construction is guaranteed by
commuta\-tivity of certain diagrams. Our main goal is to check
this directly.


We consider the two steps that seem to be necessary in the proof
of correctness of the reduced bar construction. These steps are
roughly sketched below and precisely given in Sections 4, 6 and~8.
It turns out that the definition of $n$-fold monoidal category
given in \cite{AM10} provides these two steps. We start with
checking the correctness of the reduced bar construction based on
a two-fold monoidal category, i.e.,\ 2-monoidal category of
\cite{AM10}, or duoidal category of \cite{BM12} and \cite{BS13}.
The first step in this case is trivial, and the second step, which
may be simply modified and used for the $n$-fold case, is more
involved.

Then we check the correctness of the reduced bar construction
based on a three-fold monoidal category, i.e.,\ 3-monoidal
category of \cite{AM10}. We go through two steps that are in
spirit the same as in the two-fold case. Neither of these steps is
now trivial but, as mentioned above, the second is just a
modification of the corresponding step in the two-fold case. The
combinatorial structure of $n$-fold monoidal categories, defined
by iterating this procedure, as it is already shown in
\cite{AM10}, stabilizes at $n=3$. Hence, an $n$-fold monoidal
category, for $n\geq 3$, may be envisaged as a sequence of $n$
monoidal structures in a category, such that every triple of these
structures corresponds to a three-fold monoidal category. The
correctness of the reduced bar construction based on an $n$-fold
monoidal category is obtained as a simple modification of the
results mentioned above.

\vspace{.3ex}

Our techniques are very much syntactical. We rely on a syntactical
nature of the simplicial category presented by its generating
arrows and equations. These equations are easily turned into
rewrite rules, which are useful for some normalization techniques.
Also, we try to point out the combinatorial core of the subject.
This is the reason why our definition of the reduced bar
construction $\WM$, although it covers the one of \cite{BFSV}, is
given in different terms. From a composition of functors involved
in the definition of $\WM$ we abstract a shuffle of $n$ sequences,
whose members are generators of the simplicial category. Then we
consider some available transpositions turning this shuffle into
one obtained by concatenating these $n$ sequences in a desired
order. The first step in the proof of correctness of the reduced
bar construction shows that the equations of $n$-fold monoidal
categories suffice to consider any two applications of available
transpositions from one shuffle to the other to be equal. This is
a consequence of some naturality assumptions in the two-fold case.
In the $n$-fold case, for $n\geq 3$, we need some additional
equations brought by the assumptions on 1-cells of $Mon_2(Cat)$.
Roughly speaking, these equations guarantee that the following two
applications of transpositions in our shuffles, which correspond
to the Yang-Baxter equation, are equal.

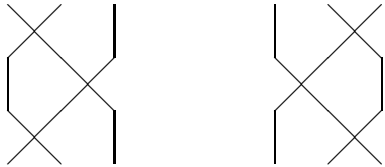
\begin{center}
\begin{picture}(140,60)

\put(0,0){\line(1,1){20}} \put(20,0){\line(-1,1){20}}
\put(40,0){\line(0,1){20}}

\put(20,20){\line(1,1){20}} \put(40,20){\line(-1,1){20}}
\put(0,20){\line(0,1){20}}

\put(0,40){\line(1,1){20}} \put(20,40){\line(-1,1){20}}
\put(40,40){\line(0,1){20}}

\put(100,20){\line(1,1){20}} \put(120,20){\line(-1,1){20}}
\put(140,20){\line(0,1){20}}

\put(120,0){\line(1,1){20}} \put(140,0){\line(-1,1){20}}
\put(100,0){\line(0,1){20}}

\put(120,40){\line(1,1){20}} \put(140,40){\line(-1,1){20}}
\put(100,40){\line(0,1){20}}

\end{picture}
\end{center}

The sequences that constitute a shuffle may be transformed
according to the equations of the simplicial category. Let $\Phi'$
be the result of such a transformation of a sequence $\Phi$. The
second step in the proof of correctness of the reduced bar
construction shows that the equations of $n$-fold monoidal
categories suffice to consider the permutation of $\Phi$ or of
$\Phi'$ with a member of another sequence to be equal. All these
equations are already present in the two-fold case.

Hence, the equations of $n$-fold monoidal categories guarantee the
correctness of the reduced bar construction. On the other hand,
these equations are also necessary if one proves the correctness
through these two steps. Our work may be characterized as the
process of defining the $n$-fold monoidal categories just from the
correctness of the reduced bar construction based on a multiple
monoidal structure. We believe there are no further possibilities
to laxify the notion of an $n$-fold monoidal category preserving
the correctness of the reduced bar construction.

With respect to the reduced bar construction, our result
generalizes all the results mentioned above. It does not involve
coherence results whose proofs in the case of $n$-fold monoidal
categories are lengthy and complicated. The two steps of our
proofs mentioned above are pretty straightforward. This paper,
except for some basic categorial definitions and results needed
for Section~9, is self-contained.

\vspace{.3ex}

To conclude, we mention that the interchanges between the monoidal
structures required for $n$-fold monoidal categories are usually
brought about by braiding and symmetry. It is pointed out in
\cite{DPInt} and \cite{PT13} that a bicartesian structure (a
category with all finite coproducts and products) brings the
desired interchanges but the corresponding coherence result
required some unusual properties of such a category---a coproduct
of terminal objects should be terminal. Our results show that this
coherence is not necessary anymore and that every bicartesian
category, for every $n$, may be conceived as an $n$-fold monoidal
category in $n+1$ different ways. Although a bicartesian category
is already $\infty$-monoidal, since it is symmetric monoidal (in
two ways), this fact is interesting---there is a family, indexed
by pairs of natural numbers, of reduced bar constructions based on
such a category. We discuss these matters in more details at the
end of Section~9. Also, this gives a positive answer to the second
question of \cite[Section~8]{PT13}.

\section{The two-fold monoidal categories}

The notion of two-fold monoidal category that we use in this paper
is defined in \cite[Section~4]{G09}. It appears in \cite[Section
6.1]{AM10} under the name 2-\textit{monoidal category} and in both
\cite[Section 2.2]{BM12} and \cite[Section~3]{BS13} under the name
\textit{duoidal category}. The notion appears as the second
iterate of the inductive definition mentioned in the introduction.
It slightly generalizes the notion of \textit{bimonoidal
intermuting category} introduced in \cite[Section 12]{DPInt}. The
difference between these two notions is that, in bimonoidal
intermuting categories, the arrows $\beta$ and~$\tau$ from below
are required to be isomorphisms. The motivation behind this
invertibility requirement is a coherence result in the style of
Kelly and Mac Lane (see \cite{KML71}), which is proved in
\cite{DPInt}.

Let $Mon(Cat)$ be the 2-category whose 0-cells are the monoidal
categories, 1-cells are the monoidal functors, and 2-cells are the
monoidal transformations (see \cite[XI.2]{ML71}). The monoidal
structure of $Mon(Cat)$ is given by 2-products (see
\cite[7.4]{B94}).

\vspace{1ex}

\noindent{\sc Definition.}\quad A \emph{two-fold monoidal}
category is a pseudomonoid in $Mon(Cat)$.

\vspace{1ex}

The unfolded form of this definition is given in Section~10
(Appendix). In this paper we are interested in strict monoidal
structures and we now give a more symmetric definition of two-fold
strict monoidal categories. A \emph{two-fold strict monoidal}
category is a category $\cal M$ equipped with two strict monoidal
structures $\langle{\cal M},\otm_1, I_1\rangle$ and $\langle{\cal
M},\otm_2, I_2\rangle$ together with the arrows $\kappa\!: I_1\str
I_2$, $\beta\!: I_1\str I_1\otm_2 I_1$, $\tau\!: I_2\otm_1 I_2\str
I_2$, and a natural transformation $\iota$ given by the family of
arrows
\[
\iota_{A,B,C,D}\!:(A\otm_2 B)\otm_1(C\otm_2 D)\str(A\otm_1
C)\otm_2(B\otm_1 D),
\]
such that the following twelve equations hold:

\begin{tabbing}
\hspace{1.5em}\=(1)\hspace{1em}\=$\iota\circ(\mj\otm_1\iota)=\iota\circ(\iota\otm_1\mj)$,
\hspace{3em}\=(7)\hspace{1em}\=$(\mj\otm_2\iota)\circ\iota=(\iota\otm_2\mj)\circ\iota$,
\\[1ex]
\>(2)\>$\iota\circ(\mj\otm_1\beta)=\mj$,\>(8)\>$(\mj\otm_2\tau)\circ\iota=\mj$,
\\[1ex]
\>(3)\>$\iota\circ(\beta\otm_1\mj)=\mj$,\>(9)\>$(\tau\otm_2\mj)\circ\iota=\mj$,
\\[1ex]
\>(4)\>$\tau\circ(\mj\otm_1\tau)=\tau\circ(\tau\otm_1\mj)$,\>(10)\>$(\mj\otm_2\beta)\circ\beta=(\beta\otm_2\mj)\circ\beta$,
\\[1ex]
\>(5)\>$\tau\circ(\mj\otm_1\kappa)=\mj$,\>(11)\>$(\mj\otm_2\kappa)\circ\beta=\mj$,
\\[1ex]
\>(6)\>$\tau\circ(\kappa\otm_1\mj)=\mj$,\>(12)\>$(\kappa\otm_2\mj)\circ\beta=\mj$.
\end{tabbing}

The two-fold monoidal categories defined in
\cite[Definition~1.7]{BFSV} are the two-fold strict monoidal
categories from above in which, moreover, it is assumed that
$I_1=I_2=0$ and $\kappa=\beta=\tau=\mj_0$. (The tensors
$\otimes_1$ and $\otimes _2$ are denoted in \cite{BFSV} by
$\square_1$ and $\square_2$, while the natural transformation
$\iota$ is denoted by $\eta$.) Hence, from our list of twelve
equations, the equations (4), (5), (6), (10), (11) and (12) are
trivial, (1) is the \emph{internal associativity condition}, (7)
is the \emph{external associativity condition}, (8) and (9) make
together the \emph{internal unit condition} and (2) and (3) make
together the \emph{external unit condition} (see
\cite[Definition~1.7]{BFSV}).

Also, every braided monoidal category is a two-fold monoidal
category in our sense. Both monoidal structures of such a two-fold
monoidal category are the same, and all the $\iota$ arrows are
obtained by braiding.


\section{The reduced bar construction}

Here we will only give a definition of the reduced bar
construction based on a strict monoidal category. We refer to
\cite[Section~6]{PT13} for the complete analysis of this
construction.

Let $\Delta$ (denoted by $\Delta^+$ in \cite{ML71}) be the
topologist's \emph{simplicial} category defined as in
\cite[VII.5]{ML71} for whose arrows we take over the notation used
in that book. In order to use geometric dimension, the objects of
$\Delta$, which are the nonempty ordinals $\{1,2,3,\ldots\}$ are
rewritten as $\{0,1,2,\ldots\}$. Hence, for $n\geq 1$ and $0\leq
i\leq n$, the source of $\delta_i^n$ is $n-1$ and the target is
$n$, while for $n\geq 1$ and $0\leq i\leq n-1$, the source of
$\sigma_i^n$ is $n$ and the target is $n-1$. When we speak of
$\Delta^{\!op}$, then we denote its arrows
$(\delta_i^n)^{op}\!:n\str n-1$ by $d_i^n$ and
$(\sigma_i^n)^{op}\!:n-1\str n$ by $s_i^n$.

The arrows of $\Delta^{\!op}$ satisfy the following \emph{basic
equations}:
\[
\begin{array}{ll}
d^{n-1}_j\circ d^n_l=d^{n-1}_{l-1}\circ d^n_j, &
\mbox{when}\hspace{1em} l-1\geq j,
\\[1ex]
s^{n+1}_j\circ s^n_l=s^{n+1}_{l+1}\circ s^n_j, &
\mbox{when}\hspace{1em} l+1>j,
\end{array}
\]
\[
d^n_j\circ s^n_l=\left\{
\begin{array}{rl}
s^{n-1}_{l-1}\circ d^{n-1}_j, & \mbox{when}\hspace{1em} j\leq l-1,
\\[1ex]
\mj, & \mbox{when}\hspace{1em} l\in\{j,j-1\},
\\[1ex]
s^{n-1}_l\circ d^{n-1}_{j-1}, & \mbox{when}\hspace{1em} j\geq l+2.
\end{array}
\right .
\]

These particular equations whose left-hand sides are treated as
redexes and the right-hand sides as the corresponding contracta
serve to define the normal form (see below). The definition of the
natural transformation $\omega$ (the ultimate ingredient in our
construction) is completely based on this normal form. We use some
syntactical techniques in this paper---it is therefore important
how we represent the arrows by terms. However, we will never write
brackets to denote the association of the binary operation of
composition, and appropriate identity arrows could be considered
present in a term or deleted from it, if necessary. The following
proposition is analogous to \cite[VII.5, Proposition~2]{ML71}.

\prop{Proposition 3.1}{The category $\Delta^{\!op}$ is generated
by the arrows $d_i^n\!:n\str n-1$ for $n\geq 1$, $0\leq i\leq n$,
and $s_i^n\!:n-1\str n$ for $n\geq 1$, $0\leq i\leq n-1$, subject
to the basic equations of $\Delta^{\!op}$.}

\dkz As in the lemma preceding \cite[VII.5, Proposition~2]{ML71},
one can prove that every arrow of $\Delta^{\!op}$ has a unique
representation of the form $\mj$ or
\[
s_{l_1}\circ\ldots\circ s_{l_k}\circ d_{j_1}\circ\ldots\circ
d_{j_m},
\]
(with the superscripts omitted) for $k+m\geq 1$, $l_1>\ldots>l_k$,
$j_1\geq\ldots\geq j_m$. The basic equations of $\Delta^{\!op}$
(read from the left to the right as reduction rules) suffice to
put any composite of $d$'s and $s$'s into the above form (cf.\ the
proof of $S4_{\Box}$ Coherence in \cite[Section~3]{DP08b}). \qed

We call the arrows $\mj_n$, $d_i^n$, and $s_i^n$ \emph{basic
arrows} of $\Delta^{\!op}$. Also, we call the above representation
of an arrow $f$ of $\Delta^{\!op}$ the \emph{normal form} of $f$
and we denote it by $f^{\rm nf}$. This normal form does not
completely correspond to the normal form given in the mentioned
lema of \cite[VII.5]{ML71}---by varying the directions of the
reduction rules corresponding to the first two basic equations of
$\Delta^{\!op}$ one may obtain other possible normal forms.

\prop{Remark 3.2}{If $f_1,\ldots,f_k$ are basic, non-identity
arrows of $\Delta^{\!op}$ such that $f_k\circ\ldots\circ f_1$ is
defined and not a normal form, then there is $1\leq i\leq k-1$
such that $f_{i+1}\circ f_i$ is the left hand side of one of the
basic equations of  $\Delta^{\!op}$.}

By \cite[XI.3, Theorem~1]{ML71}, we may regard \emph{Cat} as a
strict monoidal category whose monoidal structure is given by
finite products. Let $\cM$ be a strict monoidal category, hence a
monoid in \emph{Cat}. The \emph{reduced bar construction} (see
\cite{T79}) based on $\cM$ is a functor
$\WM\!:\Delta^{\!op}\str\mbox{\emph{Cat}}$ defined as follows.
\begin{tabbing}
\hspace{1.5em}\=$\WM(n)=\cM^n$,
\\[3ex]
\hspace{1.5em}\>$\WM(d_0^n)(A_1,A_2,\ldots,A_n)=(A_2,\ldots,A_n)$,
\\[2ex]
\>$\WM(d_n^n)(A_1,\ldots,A_{n-1},A_n)=(A_1,\ldots,A_{n-1})$,
\\[2ex]
and for $1\leq i\leq n\mn 1$ and $0\leq j\leq n\mn 1$,
\\[1ex]
\>$\WM(d_i^n)(A_1,\ldots,A_i,A_{i+1},\ldots,A_n)=(A_1,\ldots,A_i\otimes
A_{i+1},\ldots,A_n)$,
\\[2ex]
\>$\WM(s_j^n)(A_1,\ldots,A_j,A_{j+1},\ldots,A_{n-1})=
(A_1,\ldots,A_j,I,A_{j+1},\ldots,A_{n-1})$,
\\[1ex]
where $\otimes$ is the tensor and $I$ is the unit of the strict
monoidal category $\cM$.
\end{tabbing}

We denote by $\WM^m\!:\Delta^{\!op}\str\mbox{\emph{Cat}}$ the
reduced bar construction based on the $m$th power of the strict
monoidal category $\cM$ (which is again a strict monoidal category
with the structure defined component-wise). When $\cal M$ is a
two-fold strict monoidal category (or an $n$-fold, in general),
then we denote by $\WM_i\!:\Delta^{\!op}\str\mbox{\emph{Cat}}$ the
reduced bar construction based on the $i$th monoidal structure of
$\cal M$. By combining these two notations,
$\WM^m_i\!:\Delta^{\!op}\str\mbox{\emph{Cat}}$ denotes the reduced
bar construction based on the $m$th power of the strict monoidal
category whose monoidal structure is the $i$th monoidal structure
of $\cM$.

\section{The two-fold reduced bar construction}

We start with a definition of the two-fold reduced bar
construction based on a two-fold strict monoidal category. This
construction corresponds to the one given in the proof of
\cite[Theorem~2.1]{BFSV}, save that the latter construction is
based on a category that is two-fold monoidal in the sense of that
paper. Then we switch to an equivalent notion, which is of a
combinatorial flavour. Such an approach is more suitable for our
techniques, and it highlights the combinatorial core of the
subject.

Let $\cal M$ be a two-fold strict monoidal category. By relying on
the structure of $\cal M$, we define a function from objects of
$(\Delta^{\!op})^2$ to objects of \emph{Cat} and a function from
arrows of $(\Delta^{\!op})^2$ to arrows of \emph{Cat}. These two
functions are both denoted by $\WM$.

\vspace{2ex}

\noindent {\sc Definition.}\quad The \emph{two-fold reduced bar
construction} $\WM$ is defined on objects of $(\Delta^{\!op})^2$
as:
\[
\WM(n,m)=\cM^{n\cdot m},
\]
and it is defined on arrows of $(\Delta^{\!op})^2$ in the
following manner.

For $f$ an arrow of $(\Delta^{\!op})$, we have
\[
\WM(f,\mj_m)=\WM_1^m(f),
\]
where $\WM_1^m$ is, according to the notation introduced at the
end of Section~3, the reduced bar construction based on $\cM^m$
with monoidal structure given by $\otimes_1$ and $I_1$. For
example, $\WM(d^3_1,\mj_2):{\cal M}^6\str {\cal M}^4$ is such that
\[
\WM(d^3_1,\mj_2)(A,B,C,D,E,F)=(A\otimes_1 C,B\otimes_1 D,E,F),
\]
while $\WM(s^3_0,\mj_2):{\cal M}^4\str {\cal M}^6$ is such that
\[
\WM(s^2_0,\mj_2)(A,B,C,D)=(I_1,I_1,A,B,C,D).
\]

For $g$ an arrow of $(\Delta^{\!op})$, we have
\[
\WM(\mj_n,g)=(\WM_2(g))^n,
\]
where $\WM_2$ is, according to the notation introduced at the end
of Section~3, the reduced bar construction based on the strict
monoidal structure given by $\otimes_2$ and $I_2$ of $\cM$. For
example, $\WM(\mj_2,d^3_1):{\cal M}^6\str {\cal M}^4$ is such that
\[
\WM(\mj_2,d^3_1)(A,B,C,D,E,F)=(A\otimes_2 B,C,D\otimes_2E,F),
\]
while $\WM(\mj_2,s^3_0):{\cal M}^4\str {\cal M}^6$ is such that
\[
\WM(\mj_2,s^3_0)(A,B,C,D)=(I_2,A,B,I_2,C,D).
\]

Finally, for $f:n_s\str n_t$ and $g:m_s\str m_t$, (``$s$'' comes
from \emph{source} and ``$t$'' from \emph{target}) we have
\[
\WM(f,g)=(\WM_2(g))^{n_t}\circ \WM_1^{m_s}(f).
\]
For example, $\WM(d^3_1,s^3_0):{\cal M}^6\str {\cal M}^6$ is such
that
\[
\WM(d^3_1,s^3_0)(A,B,C,D,E,F)=(I_2,A\otimes_1 C,B\otimes_1
D,I_2,E,F)
\]

In general, $\WM$ is not a functor from $(\Delta^{\!op})^2$ to
\emph{Cat} since it does not preserve composition (it preserves
identities). For example,
$\WM(d^2_1,\mj_1)\circ\WM(\mj_2,d^2_1):{\cal M}^4\str {\cal M}$ is
such that
\[
(\WM(d^2_1,\mj_1)\circ\WM(\mj_2,d^2_1))(A,B,C,D)=(A\otimes_2
B)\otimes_1(C\otimes_2 D),
\]
while $\WM(d^2_1,d^2_1):{\cal M}^4\str {\cal M}$ is such that
\[
\WM(d^2_1,d^2_1)(A,B,C,D)=(A\otimes_1 C)\otimes_2(B\otimes_1 D).
\]

Our goal is to show that $\WM:(\Delta^{\!op})^2\str
\mbox{\emph{Cat}}$ is a \emph{lax functor} in the sense of
\cite{S72}. This means that for every composable pair of arrows
$e_1=(f_1,g_1)$ and $e_2=(f_2,g_2)$ of $(\Delta^{\!op})^2$, there
is a natural transformation
\[
\omega_{e_2,e_1}\!:\WM(e_2)\circ \WM(e_1)\strt \WM(e_2\circ e_1),
\]
such that the following diagram commutes:

\begin{center}
\begin{picture}(100,110)(0,-20)

\put(50,75){\makebox(0,0){$\WM(e_3)\circ\WM(e_2)\circ \WM(e_1)$}}
\put(-40,40){\makebox(0,0){$\WM(e_3\circ e_2)\circ \WM(e_1)$}}
\put(140,40){\makebox(0,0){$\WM(e_3)\circ\WM(e_2\circ e_1)$}}
\put(50,10){\makebox(0,0){$\WM(e_3\circ e_2\circ e_1)$}}

\put(-10,60){\vector(-2,-1){20}} \put(110,60){\vector(2,-1){20}}
\put(-30,30){\vector(2,-1){30}} \put(130,30){\vector(-2,-1){30}}

\put(-20,60){\makebox(0,0)[r]{\scriptsize
$\omega_{e_3,e_2}\WM(e_1)$}}
\put(120,60){\makebox(0,0)[l]{\scriptsize
$\WM(e_3)\omega_{e_2,e_1}$}}

\put(-20,20){\makebox(0,0)[r]{\scriptsize $\omega_{e_3\circ
e_2,e_1}$}} \put(120,20){\makebox(0,0)[l]{\scriptsize
$\omega_{e_3,e_2\circ e_1}$}}

\put(50,-15){\makebox(0,0){(Diag 4.1)}}

\end{picture}
\end{center}
The rest of this section is devoted to a proof of laxness of
$\WM$.

\vspace{2ex}

For $k\geq 1$, let $f_k\ldots f_1$ be a sequence of basic arrows
of $\Delta^{\!op}$ such that the composition $f_k\circ\ldots\circ
f_1$ is defined. We say that $\Phi=(f_k,1)\ldots(f_1,1)$ is a
sequence of \emph{colour}~1 and we abbreviate the term
$f_k\circ\ldots\circ f_1$ by $\circ\Phi$. A sequence of
\emph{colour}~2 (or of any other colour) is defined in the same
manner. We assume that, if necessary, appropriate identities could
always be added to, or deleted from sequences of any colour.
However, for measuring the \emph{length} of such a sequence, only
non-identity members are taken into account.

Let $\Phi=(f_k,1)\ldots(f_1,1)$ be a sequence of colour~1 and let
$\Gamma=(g_l,2)\ldots(g_1,2)$ be a sequence of colour~2, such that
$\circ\Phi\!:n_s\str n_t$ and $\circ\Gamma\!:m_s\str m_t$. Let
$\cal M$ be a two-fold strict monoidal category. We define a
functor
\[
\WM_{\Gamma\Phi}\!:{\cal M}^{n_s\cdot m_s}\str{\cal M}^{n_t\cdot
m_t}
\]
as the following composition
\[
(\WM_2(g_l))^{n_t}\circ\ldots\circ(\WM_2(g_q))^{n_t}\circ
\WM_1^{m_s}(f_k)\circ\ldots\circ \WM_1^{m_s}(f_1).
\]

Let $f=\circ\Phi$ and $g=\circ\Gamma$. Since both $\WM_1$ and
$\WM_2$ are functors, we have that $\WM_{\Gamma\Phi}=\WM(f,g)$.
This fact leads to a \emph{combinatorial} definition of the
two-fold reduced bar construction $\WM$, according to which
$\WM(f,g)$ could be defined as $\WM_{\Gamma\Phi}$ for arbitrary
$\Phi$ of colour~1 such that $\circ\Phi=f$ and arbitrary $\Gamma$
of colour~2 such that $\circ\Gamma=g$.

In order to define the natural transformations $\omega$ involved
in Diag~4.1, we introduce the following notions. Let $\Theta$ be a
\textit{shuffle} of $\Phi$ and $\Gamma$ as above. For example, let
$\Phi$ be $(d^2_1,1)(d^3_1,1)$, let $\Gamma$ be
$(d^3_2,2)(s^3_0,2)(d^3_1,2)$, and let $\Theta$ be the following
shuffle
\[
(d^3_2,2)(d^2_1,1)(d^3_1,1)(s^3_0,2)(d^3_1,2).
\]
In this case, we have that $\circ\Phi\!:3\str 1$ and
$\circ\Gamma\!:3\str 2$.

For every member $(f,1)$ of~$\Theta$, we define its \emph{inner
power} to be the target of its right-closest $(g,2)$ in~$\Theta$.
We may assume that such $(g,2)$ exists since we can always add the
appropriate identity of colour 2 to the right of $(f,1)$
in~$\Theta$. For every member $(g,2)$ of $\Theta$, we define its
\emph{outer power} to be the target of its right-closest $(f,1)$
in~$\Theta$. For $\Theta$ as above, we have that the inner power
of $(d^2_1,1)$ is 3 and the outer power of $(d^3_2,2)$ is 1.

For a two-fold strict monoidal category $\cal M$ and for an
arbitrary shuffle $\Theta$ of $\Phi$ and $\Gamma$, as for the
shuffle $\Gamma\Phi$ (obtained by concatenating $\Gamma$ and
$\Phi$), we can define a functor
\[
\WM_\Theta\!:{\cal M}^{n_s\cdot m_s}\str{\cal M}^{n_t\cdot m_t}
\]
in the following way: replace in $\Theta$ every $(f,1)$ whose
inner power is $i$ by $\WM_1^i(f)$, and every $(g,2)$ whose outer
power is $o$ by $(\WM_2(g))^o$, and insert $\circ$'s. For $\Theta$
as above, we have that $\WM_\Theta$~is
\[
\WM_2(d^3_2)\circ\WM_1^3(d^2_1)\circ\WM_1^3(d^3_1)
\circ(\WM_2(s^3_0))^3\circ(\WM_2(d^3_1))^3,
\]
which gives that $\WM_\Theta(A,B,C,D,E,F,G,H,J)$ is the ordered
pair
\[
(I_2\otm_1 I_2\otm_1 I_2,((A\otm_2 B)\otm_1(D\otm_2
E)\otm_1(G\otm_2 H))\otm_2(C\otm_1 F\otm_1 J)).
\]

For basic arrows $f\!:n\str n'$ and $g\!:m\str m'$ of
$\Delta^{\!op}$ we define a natural transformation
\[
\chi(f,g)\!:\WM_1^{m'}(f)\circ(\WM_2(g))^n\;\strt\;
(\WM_2(g))^{n'}\circ\WM_1^m(f)
\]
to be the identity natural transformation except in the following
cases:

\begin{center}

\begin{tabular}{c|c|c}
$f$ & $g$ & $\chi(f,g)$
\\[.5ex]
\hline\hline\\[-2.2ex] $s^{n+1}_j$ & $s^{m+1}_i$ &
$(\mj^{j(m+1)},\mj^i,\kappa,\vmj)$
\\[.5ex]
\hline\\[-2.2ex] $d^n_j$, $1\leq j\leq n-1$ & $s^{m+1}_i$ &
$(\mj^{(j-1)(m+1)},\mj^i,\tau,\vmj)$
\\[.5ex]
\hline\\[-2.2ex] $s^{n+1}_j$ & $d^m_i$, $1\leq i\leq m-1$ &
$(\mj^{j(m-1)},\mj^{i-1},\beta,\vmj)$
\\[.5ex]
\hline\\[-2.2ex] $d^n_j$, $1\leq j\leq n-1$ & $d^m_i$, $1\leq i\leq m-1$ &
$(\mj^{(j-1)(m-1)},\mj^{i-1},\iota,\vmj)$
\end{tabular}
\captionof{table}{$\chi(f,g)$ in nontrivial cases.}
\label{table:2dim}
\end{center}
Here $\mj^n$ denotes the $n$-tuple of identities and $\vmj$ is a
tuple of identities whose length can be easily calculated in all
the cases, but we will not write the exact length to avoid
overlong expressions.

For $j\geq 0$, let $\Theta_0,\ldots,\Theta_j$ be shuffles of
$\Phi$ and $\Gamma$ such that $\Theta_0$ is $\Theta$ and
$\Theta_j$ is $\Gamma\Phi$, and if $j>0$, then for every $0\leq
i\leq j-1$ we have that for some shuffles $\Pi$ and $\Lambda$ and
non-identity members $(f,1),(g,2)$, the shuffle $\Theta_i$ is
$\Pi(f,1)(g,2)\Lambda$ while $\Theta_{i+1}$ is
$\Pi(g,2)(f,1)\Lambda$. We call $\Theta_0,\ldots,\Theta_j$ a
\emph{normalizing path} starting with $\Theta$. Its \emph{length}
is $j$. For example,
\begin{tabbing}
\hspace{1em}\=$\Theta_0=(d^3_2,2)(d^2_1,1)(d^3_1,1)(s^3_0,2)(d^3_1,2)$,
$\Theta_1=(d^3_2,2)(d^2_1,1)(s^3_0,2)(d^3_1,1)(d^3_1,2)$,
\\[1ex]
\>$\Theta_2=(d^3_2,2)(d^2_1,1)(s^3_0,2)(d^3_1,2)(d^3_1,1)$,
$\Theta_3=(d^3_2,2)(s^3_0,2)(d^2_1,1)(d^3_1,2)(d^3_1,1)$,
\\[1ex]
\>$\Theta_4=(d^3_2,2)(s^3_0,2)(d^3_1,2)(d^2_1,1)(d^3_1,1)$
\end{tabbing}
is a normalizing path of length 4 starting with $\Theta$ as in the
example given above.

\prop{Proposition 4.1}{Every normalizing path starting with
$\Theta$ has the same length.}

\dkz The length of every normalizing path starting with $\Theta$
is $\displaystyle \sum_{(f,1)\:{\rm\scriptstyle
in}\:\Theta}k(f,1)$, where $k(f,1)$ is the number of non-identity
members of colour~2 to the right of $(f,1)$ in $\Theta$. \qed

If $\Theta_i=\Pi(f,1)(g,2)\Lambda$ and
$\Theta_{i+1}=\Pi(g,2)(f,1)\Lambda$, then
\[
\varphi_i=\WM_\Pi\;\chi(f,g)\;\WM_\Lambda
\]
is a natural transformation from $\WM_{\Theta_i}$ to
$\WM_{\Theta_{i+1}}$. (In the case when $\Pi$ or $\Lambda$ are
single-coloured, we can always add the appropriate identity of the
other colour in order to define $\WM_\Pi$ and $\WM_\Lambda$.) Let
\[
\varphi(\Theta_0,\ldots,\Theta_j)=\left\{
\begin{array}{rl}
\varphi_{j-1}\circ\ldots\circ\varphi_0, & \mbox{when}\hspace{1em}
j\geq 1,
\\[1ex]
\mj, & \mbox{when}\hspace{1em} j=0.
\end{array}
\right .
\]

Suppose $\Theta'_0,\ldots,\Theta'_j$ is another normalizing path
starting with $\Theta$. Then $\varphi(\Theta'_0,\ldots,\Theta'_j)$
is again a natural transformation from $\WM_\Theta$ to
$\WM_{\Gamma\Phi}$. We can show that these natural transformations
are in fact the same.

\prop{Theorem 4.2}{
$\varphi(\Theta_0,\ldots,\Theta_j)=\varphi(\Theta'_0,\ldots,\Theta'_j)$.}

\dkz We proceed by induction on $j\geq 0$. If $j=0$, then
$\varphi(\Theta_0)=\varphi(\Theta'_0)=\mj$. If $j>0$ and
$\Theta_1=\Theta_1'$, then we apply the induction hypothesis to
the sequences of shuffles $\Theta_1,\ldots,\Theta_j$ and
$\Theta'_1,\ldots,\Theta'_j$. Suppose now  $\Theta_1\neq\Theta_1'$
and $\varphi_0=\WM_\Pi\;\chi(f,g)\;\WM_\Lambda$,
$\varphi'_0=\WM_{\Pi'}\;\chi(f',g')\;\WM_{\Lambda'}$. Then either
\[
\Theta=\Pi_1(f',1)(g',2)\Pi_2(f,1)(g,2)\Lambda\quad {\rm or}\quad
\Theta=\Pi(f,1)(g,2)\Lambda_1(f',1)(g',2)\Lambda_2.
\]
In the first case, by naturality we have
\begin{tabbing}
\hspace{.5em}$(\ast)$\quad
$\WM_{\Pi_1}\;\chi(f',g')\;\WM_{\Pi_2(g,2)(f,1)\Lambda}\circ
\WM_{\Pi_1(f',1)(g',2)\Pi_2}\;\chi(f,g)\;\WM_\Lambda$
\\[1ex]
\`$= \WM_{\Pi_1(g',2)(f',1)\Pi_2}\;\chi(f,g)\;\WM_\Lambda\circ
\WM_{\Pi_1}\;\chi(f',g')\;\WM_{\Pi_2(f,1)(g,2)\Lambda}$,
\end{tabbing}
and by applying the induction hypothesis twice we obtain the
following commutative diagram, in which $\Xi$ is
$\Pi_1(g',2)(f',1)\Pi_2(g,2)(f,1)\Lambda$.

\begin{center}
\begin{picture}(100,140)(0,-5)

\put(55,-5){\makebox(0,0){$\WM_{\Gamma\Phi}$}}
\put(53,55){\makebox(0,0){$\WM_{\Xi}$}}
\put(0,90){\makebox(0,0){$\WM_{\Theta_1}$}}
\put(105,90){\makebox(0,0){$\WM_{\Theta_1'}$}}
\put(50,128){\makebox(0,0){$\WM_{\Theta}$}}

\put(20,115){\makebox(0,0)[r]{$\varphi_0$}}
\put(80,115){\makebox(0,0)[l]{$\varphi_0'$}}
\put(-3,70){\makebox(0,0)[r]{$\varphi_1$}}
\put(103,70){\makebox(0,0)[l]{$\varphi_1'$}}
\put(20,10){\makebox(0,0)[r]{$\varphi_{j-1}$}}
\put(80,10){\makebox(0,0)[l]{$\varphi_{j-1}'$}}

\put(50,90){\makebox(0,0){$(\ast)$}}
\put(25,40){\makebox(0,0){\scriptsize ind. hyp.}}
\put(75,40){\makebox(0,0){\scriptsize ind. hyp.}}

\put(0,45){\makebox(0,0){$\vdots$}}
\put(50,45){\makebox(0,0){$\vdots$}}
\put(100,45){\makebox(0,0){$\vdots$}}

\put(38,118){\vector(-3,-2){30}} \put(62,118){\vector(3,-2){30}}
\put(92,83){\vector(-3,-2){30}} \put(8,83){\vector(3,-2){30}}
\put(92,23){\vector(-3,-2){30}} \put(8,23){\vector(3,-2){30}}

\put(50,30){\vector(0,-1){20}} \put(0,80){\vector(0,-1){20}}
\put(100,80){\vector(0,-1){20}}
\end{picture}
\end{center}
We proceed analogously in the second case. \qed

By Theorem 4.2, the following definition is correct.

\vspace{2ex}

\noindent {\sc Definition.}\quad Let $\varphi_\Theta:\WM_\Theta
\strt \WM_{\Gamma\Phi}$ be $\varphi(\Theta_0,\ldots,\Theta_j)$,
for an arbitrary normalizing path $\Theta_0,\ldots,\Theta_j$
starting with $\Theta$.

\vspace{2ex}

For every composable pair of arrows $e_1=(f_1,g_1)$ and
$e_2=(f_2,g_2)$ of $(\Delta^{\!op})^2$ we define a natural
transformation
\[
\omega_{e_2,e_1}\!:\WM(e_2)\circ \WM(e_1)\strt \WM(e_2\circ e_1).
\]
In order to do this, note that for a sequence H of any colour,
$\circ{\rm H}$ denotes a syntactical object, a word of the form
$h_k\circ\ldots\circ h_1$. Hence, a sequence H is completely
determined by its colour and $\circ{\rm H}$.

\vspace{2ex}

\noindent {\sc Definition.}\quad Let $\Phi_1$ and $\Phi_2$ be
sequences of colour 1, and let $\Gamma_1$ and $\Gamma_2$ be
sequences of colour 2 such that $\circ\Phi_1$ is $f^{\rm nf}_1$,
$\circ\Phi_2$ is $f^{\rm nf}_2$, $\circ\Gamma_1$ is $g^{\rm nf}_1$
and $\circ\Gamma_2$ is $g^{\rm nf}_2$. We define
\[
\omega_{e_2,e_1}\quad {\rm as}\quad
\varphi_{\Gamma_2\Phi_2\Gamma_1\Phi_1}.
\]

\noindent\textit{Note.} The source and target of
$\omega_{e_2,e_1}$ are as desired since
\[
\WM(e_2)\circ\WM(e_1)=\WM_{\Gamma_2\Phi_2\Gamma_1\Phi_1},\;
\mbox{\rm and}
\]
\[
\WM(e_2\circ e_1)=\WM_{\Gamma_2\Gamma_1\Phi_2\Phi_1}.
\]

It remains to prove that our Diag~4.1 commutes. Let
$e_1=(f_1,g_1)$, $e_2=(f_2,g_2)$ and $e_3=(f_3,g_3)$ be such that
the composition $e_3\circ e_2\circ e_1$ is defined in
$(\Delta^{\!op})^2$. Let $\Phi_1$, $\Phi_2$, $\Gamma_1$ and
$\Gamma_2$ be as above, and let $\Phi_3$ and $\Gamma_3$ be
sequences of colour~1 and 2 respectively such that $\circ\Phi_3$
is~$f^{\rm nf}_3$ and $\circ\Gamma_3$ is~$g^{\rm nf}_3$. In this
case, Diag~4.1 reads

\begin{center}
\begin{picture}(100,70)(0,5)

\put(50,75){\makebox(0,0){$\WM(e_3)\circ\WM(e_2)\circ \WM(e_1)$}}
\put(-40,40){\makebox(0,0){$\WM(e_3\circ e_2)\circ \WM(e_1)$}}
\put(140,40){\makebox(0,0){$\WM(e_3)\circ\WM(e_2\circ e_1)$}}
\put(50,10){\makebox(0,0){$\WM(e_3\circ e_2\circ e_1)$}}

\put(-10,60){\vector(-2,-1){20}} \put(110,60){\vector(2,-1){20}}
\put(-30,30){\vector(2,-1){30}} \put(130,30){\vector(-2,-1){30}}

\put(-20,60){\makebox(0,0)[r]{\scriptsize
$\varphi_{\Gamma_3\Phi_3\Gamma_2\Phi_2}\WM_{\Gamma_1\Phi_1}$}}
\put(120,60){\makebox(0,0)[l]{\scriptsize $\WM_{\Gamma_3\Phi_3}\;
\varphi_{\Gamma_2\Phi_2\Gamma_1\Phi_1}$}}

\put(-20,20){\makebox(0,0)[r]{\scriptsize
$\varphi_{\Gamma'\Phi'\Gamma_1\Phi_1}$}}
\put(120,20){\makebox(0,0)[l]{\scriptsize
$\varphi_{\Gamma_3\Phi_3\Gamma''\Phi''}$}}

\end{picture}
\end{center}
where $\circ\Phi'$ is $(f_3\circ f_2)^{\rm nf}$, $\circ\Gamma'$ is
$(g_3\circ g_2)^{\rm nf}$, $\circ\Phi''$ is $(f_2\circ f_1)^{\rm
nf}$ and $\circ\Gamma''$ is $(g_2\circ g_1)^{\rm nf}$.

\vspace{2ex}

By Theorem 4.2 we have the following commutative diagram

\begin{center}
\begin{picture}(100,75)(0,5)

\put(50,75){\makebox(0,0){$\WM(e_3)\circ\WM(e_2)\circ \WM(e_1)$}}
\put(-40,40){\makebox(0,0){$\WM(e_3\circ e_2)\circ \WM(e_1)$}}
\put(140,40){\makebox(0,0){$\WM(e_3)\circ\WM(e_2\circ e_1)$}}
\put(50,10){\makebox(0,0){$\WM(e_3\circ e_2\circ e_1)$}}

\put(-10,60){\vector(-2,-1){20}} \put(110,60){\vector(2,-1){20}}
\put(-30,30){\vector(2,-1){30}} \put(130,30){\vector(-2,-1){30}}
\put(50,65){\line(0,-1){7}} \put(50,50){\vector(0,-1){30}}

\put(-20,60){\makebox(0,0)[r]{\scriptsize
$\varphi_{\Gamma_3\Phi_3\Gamma_2\Phi_2}\WM_{\Gamma_1\Phi_1}$}}
\put(120,60){\makebox(0,0)[l]{\scriptsize $\WM_{\Gamma_3\Phi_3}\;
\varphi_{\Gamma_2\Phi_2\Gamma_1\Phi_1}$}}

\put(-20,20){\makebox(0,0)[r]{\scriptsize
$\varphi_{\Gamma_3\Gamma_2\Phi_3\Phi_2\Gamma_1\Phi_1}$}}
\put(120,20){\makebox(0,0)[l]{\scriptsize
$\varphi_{\Gamma_3\Phi_3\Gamma_2\Gamma_1\Phi_2\Phi_1}$}}
\put(50,55){\makebox(0,0){\scriptsize
$\varphi_{\Gamma_3\Phi_3\Gamma_2\Phi_2\Gamma_1\Phi_1}$}}

\end{picture}
\end{center}
Hence, to prove that Diag~4.1 commutes, it suffices to show that
\[
(i)\quad
\varphi_{\Gamma_3\Gamma_2\Phi_3\Phi_2\Gamma_1\Phi_1}=\varphi_{\Gamma'\Phi'\Gamma_1\Phi_1}
\quad {\rm and}\quad (ii)\quad
\varphi_{\Gamma_3\Phi_3\Gamma_2\Gamma_1\Phi_2\Phi_1}=\varphi_{\Gamma_3\Phi_3\Gamma''\Phi''}.
\]


\prop{Lemma 4.3}{If $\Phi$ and $\Phi'$ are sequences of colour 1
such that $\circ\Phi=\circ\Phi'$ is a basic equation of
$\Delta^{\!op}$, and $g$ is a basic arrow of $\Delta^{\!op}$, then
$\varphi_{\Phi (g,2)}=\varphi_{\Phi' (g,2)}$.}

\dkz We have the following cases in which we always assume that
$d^x_y$ is such that $1\leq y\leq x-1$ (see
Table~\ref{table:2dim}). To deal with $d^x_0$ and $d^x_x$ is
trivial. We will give a detailed proof for three cases, first of
which is trivial, with the remaining two needing some of the
equations (1)---(6). The rest is done analogously.

\vspace*{3pt}

\noindent 1.1. Suppose $\circ\Phi=\circ\Phi'$ is $d^{n-1}_j\circ
d^n_l=d^{n-1}_{l-1}\circ d^n_j$ for $j\leq l-2$.

\vspace{1ex}

1.1.1. Suppose $g$ is $s^m_i$.

\noindent We have two normalizing paths. The first one is starting
with ${\Phi (g,2)}$ and it is
\[
 (d_j^{n-1},1)(d_{l}^n,1)(s_i^m,2),\ (d_j^{n-1},1)(s_i^m,2)(d_{l}^n,1),\ (s_i^m,2)(d_j^{n-1},1)(d_{l}^n,1).
\]
%
Now we compute $\varphi_{\Phi (g,2)}$, and we note that $\WM_{(d_{l}^n,1)}$ is formally $\WM_{(d_{l}^n,1),(\mathbf{1}_{m-1},2)}$ (we repeatedly use such an abbreviation throughout the paper): 
\begin{equation*}
\begin{split}
 \varphi_{\Phi (g,2)} &=  \left(\chi(d_j^{n-1}, s_i^m)\;\WM_{(d_{l}^n,1)}\right) \circ  \left(\WM_{(d_{j}^{n-1},1)}\;\chi(d_{l}^{n}, s_i^m)\right)\\
   &= \left(\chi(d_j^{n-1}, s_i^m)\; \WM_1^{m-1}(d_{l}^n)\right) \circ  \left(\WM_1^m(d_{j}^{n-1})\; \chi(d_{l}^{n}, s_i^m)\right)\\
   &= \left((\mj^{(j-1)m},\mj^i,\tau,\vmj)\; \WM_1^{m-1}(d_{l}^n)\right) \circ  \left(\WM_1^m(d_{j}^{n-1})\;(\mj^{(l-1)m},\mj^i,\tau,\vmj)\right)\\
   &= (\mj^{(j-1)m},\mj^i,\tau,\vmj) \circ (\mj^{(l-2)m},\mj^i,\tau,\vmj).\qquad(\textnormal{since }l-1>j)
\end{split}
\end{equation*}
On the other hand, the second normalizing path starting with
${\Phi' (g,2)}$  is
\[ (d_{l-1}^{n-1},1)(d_{j}^n,1)(s_i^m,2), \  (d_{l-1}^{n-1},1)(s_i^m,2)(d_{j}^n,1),  \ (s_i^m,2)(d_{l-1}^{n-1},1)(d_{j}^n,1),\]
and therefore
\begin{equation*}
\begin{split}
 \varphi_{\Phi' (g,2)} &=  \left(\chi(d_{l-1}^{n-1}, s_i^m)\; \WM_{(d_{j}^n,1)}\right) \circ  \left(\WM_{(d_{l-1}^{n-1},1)}\; \chi(d_{j}^{n}, s_i^m)\right)\\
   &= \left(\chi(d_{l-1}^{n-1}, s_i^m)\; \WM_1^{m-1}(d_{j}^n)\right) \circ  \left(\WM_1^m(d_{l-1}^{n-1})\; \chi(d_{j}^{n}, s_i^m)\right)\\
   &= \left((\mj^{(l-2)m},\mj^i,\tau,\vmj)\; \WM_1^{m-1}(d_{j}^n)\right) \circ  \left(\WM_1^m(d_{l-1}^{n-1})\;(\mj^{(j-1)m},\mj^i,\tau,\vmj)\right)\\
   &= (\mj^{(l-2)m},\mj^i,\tau,\vmj) \circ (\mj^{(j-1)m},\mj^i, \tau,\vmj).
\end{split}
\end{equation*}
Since $j-1 \neq l-2$, we see that these two tuples of arrows are
the same, i.e.,\ we have:
\begin{tabbing}
\hspace{3em}\=$\varphi_{\Phi (g,2)}=(\mj^{(j-1)
m+i},\tau,\mj^{(l-j-1)m-1},\tau,\vmj)=\varphi_{\Phi' (g,2)}$.
\end{tabbing}

\vspace{1ex}

1.1.2. Suppose $g$ is $d^m_i$.
\begin{tabbing}
\hspace{3em}\=$\varphi_{\Phi (g,2)}=(\mj^{(j-1)
(m-1)+i-1},\iota,\mj^{(l-j-1)(m-1)-1},\iota,\vmj)=\varphi_{\Phi'
(g,2)}$.
\end{tabbing}

\vspace{2ex}

\noindent 1.2. Suppose $\circ\Phi=\circ\Phi'$ is $d^{n-1}_j\circ
d^n_{j+1}=d^{n-1}_{j}\circ d^n_j$.

\vspace{1ex}

1.2.1. Suppose $g$ is $s^m_i$.\\

The normalizing path starting with ${\Phi (g,2)}$ is
\[(d_j^{n-1},1)(d_{j+1}^n,1)(s_i^m,2),\
(d_j^{n-1},1)(s_i^m,2)(d_{j+1}^n,1),\
(s_i^m,2)(d_j^{n-1},1)(d_{j+1}^n,1),\] and we have
\begin{equation*}
\begin{split}
 \varphi_{\Phi (g,2)} &=  \left(\chi(d_j^{n-1}, s_i^m)\;\WM_{(d_{j+1}^n,1)}\right) \circ  \left(\WM_{(d_{j}^{n-1},1)}\;\chi(d_{j+1}^{n}, s_i^m)\right)\\
   &= \left(\chi(d_j^{n-1}, s_i^m)\; \WM_1^{m-1}(d_{j+1}^n)\right) \circ  \left(\WM_1^m(d_{j}^{n-1})\; \chi(d_{j+1}^{n}, s_i^m)\right)\\
   &= \left((\mj^{(j-1)m},\mj^i,\tau,\vmj)\; \WM_1^{m-1}(d_{j+1}^n)\right) \circ  \left(\WM_1^m(d_{j}^{n-1})\;(\mj^{jm},\mj^i,\tau,\vmj)\right)\\
   &= (\mj^{(j-1)m},\mj^i,\tau,\vmj) \circ (\mj^{(j-1)m},\mj^i,\mj\otimes_1 \tau,\vmj)\\ &= (\mj^{(j-1)m+i},\tau\circ(\mj\otimes_1\tau),\vmj).
\end{split}
\end{equation*}
On the other hand, the normalizing path starting with ${\Phi'
(g,2)}$ is
\[(d_j^{n-1},1)(d_{j}^n,1)(s_i^m,2),\ (d_j^{n-1},1)(s_i^m,2)(d_{j}^n,1),\ (s_i^m,2)(d_j^{n-1},1)(d_{j}^n,1).\] We now compute $\varphi_{\Phi' (g,2)}$:
\begin{equation*}
\begin{split}
 \varphi_{\Phi' (g,2)} &=  \left(\chi(d_j^{n-1}, s_i^m)\; \WM_{(d_{j}^n,1)}\right) \circ  \left(\WM_{(d_{j}^{n-1},1)}\; \chi(d_{j}^{n}, s_i^m)\right)\\
   &= \left(\chi(d_j^{n-1}, s_i^m)\; \WM_1^{m-1}(d_{j}^n)\right) \circ  \left(\WM_1^m(d_{j}^{n-1})\; \chi(d_{j}^{n}, s_i^m)\right)\\
   &= \left((\mj^{(j-1)m},\mj^i,\tau,\vmj)\; \WM_1^{m-1}(d_{j}^n)\right) \circ  \left(\WM_1^m(d_{j}^{n-1})\;(\mj^{(j-1)m},\mj^i,\tau,\vmj)\right)\\
   &= (\mj^{(j-1)m},\mj^i,\tau,\vmj) \circ (\mj^{(j-1)m},\mj^i, \tau \otimes_1 \mj,\vmj)\\ &= (\mj^{(j-1)m+i},\tau\circ(\tau\otimes_1\mj),\vmj).
\end{split}
\end{equation*}
Since, by (4), we have that $\tau\circ(\mj\otimes_1\tau) =
\tau\circ(\tau\otimes_1\mj)$, we conclude that $\varphi_{\Phi
(g,2)}=\varphi_{\Phi' (g,2)}$.

\vspace{1ex}

1.2.2. Suppose $g$ is $d^m_i$.
\begin{tabbing}
\hspace{3em}\=$\varphi_{\Phi (g,2)}$ \=$=(\mj^{(j-1)
(m-1)+i-1},\iota\circ(\mj\otimes_1\iota),\vmj)$
\\[.5ex]
\>\>$=(\mj^{(j-1)(m-1)+i-1},\iota\circ(\iota\otimes_1\mj),\vmj)=\varphi_{\Phi'
(g,2)}$,\quad by (1).
\end{tabbing}

\vspace{2ex}

\noindent 2. Suppose $\circ\Phi=\circ\Phi'$ is $s^{n+1}_j\circ
s^n_l=s^{n+1}_{l+1}\circ s^n_j$ for $j\leq l$.

\vspace{1ex}

2.1. Suppose $g$ is $s^m_i$.
\begin{tabbing}
\hspace{3em}\=$\varphi_{\Phi
(g,2)}=(\mj^{jm+i},\kappa,\mj^{(l-j+1)m-1},\kappa,\vmj)=\varphi_{\Phi'
(g,2)}$.
\end{tabbing}

\vspace{1ex}

2.2. Suppose $g$ is $d^m_i$.
\begin{tabbing}
\hspace{3em}\=$\varphi_{\Phi
(g,2)}=(\mj^{j(m-1)+i-1},\beta,\mj^{(l-j+1)(m-1)-1},\beta,\vmj)=\varphi_{\Phi'
(g,2)}$.
\end{tabbing}

\vspace{2ex}

\noindent 3.1. Suppose $\circ\Phi=\circ\Phi'$ is $d^n_j\circ
s^n_l=s^{n-1}_{l-1}\circ d^{n-1}_j$ for $j\leq l-1$.

\vspace{1ex}

3.1.1. Suppose $g$ is $s^m_i$.
\begin{tabbing}
\hspace{3em}\=$\varphi_{\Phi
(g,2)}=(\mj^{(j-1)m+i},\tau,\mj^{(l-j)m-1},\kappa,\vmj)=\varphi_{\Phi'
(g,2)}$.
\end{tabbing}

\vspace{1ex}

3.1.2. Suppose $g$ is $d^m_i$.
\begin{tabbing}
\hspace{3em}\=$\varphi_{\Phi
(g,2)}=(\mj^{(j-1)(m-1)+i-1},\iota,\mj^{(l-j)(m-1)-1},\beta,\vmj)=\varphi_{\Phi'
(g,2)}$.
\end{tabbing}

\vspace{2ex}

\noindent 3.2. Suppose $\circ\Phi=\circ\Phi'$ is $d^n_j\circ
s^n_j=\mj$.

\vspace{1ex}

3.2.1. Suppose $g$ is $s^m_i$.
\begin{tabbing}
\hspace{3em}\=$\varphi_{\Phi (g,2)}$ \=$=(\mj^{(j-1)
m+i},\tau\circ(\mj\otimes_1\kappa),\vmj)$
\\[.5ex]
\>\>$=(\mj^{(j-1)m+i},\mj,\vmj)=\varphi_{\Phi' (g,2)}$,\quad by
(5).
\end{tabbing}

\vspace{1ex}

3.2.2. Suppose $g$ is $d^m_i$.\\
We now have this normalizing path starting with $\Phi (g,2)$:
\[
 (d_j^{n},1)(s_{j}^n,1)(d_i^m,2),\quad (d_j^{n},1)(d_i^m,2)(s_{j}^n,1), \quad (d_i^m,2)(d_j^{n},1)(s_{j}^n,1).
\]
Since $\varphi_{\Phi' (g,2)} = \varphi_{\mathbf{1} (g,2)} = \vmj$,
we ought to compute $\varphi_{\Phi (g,2)}$:
\begin{equation*}
\begin{split}
 \varphi_{\Phi (g,2)} &=  \left(\chi(d_j^{n}, d_i^m)\; \WM_{(s_{j}^n,1)}\right) \circ  \left(\WM_{(d_{j}^{n},1)}\; \chi(s_{j}^{n}, d_i^m)\right)\\
   &= \left(\chi(d_j^{n}, d_i^m)\; \WM_1^{m}(s_{j}^n)\right) \circ  \left(\WM_1^{m-1}(d_{j}^{n})\; \chi(s_{j}^{n}, d_i^m)\right)\\
   &= \left((\mj^{(j-1)(m-1)},\mj^{i-1},\iota,\vmj) \WM_1^{m}(s_{j}^n)\right) \circ  \left(\WM_1^{m-1}(d_{j}^{n}) (\mj^{j(m-1)},\mj^{i-1},\beta,\vmj)\right)\\
   &= (\mj^{(j-1)(m-1)},\mj^{i-1},\iota,\vmj) \circ (\mj^{(j-1)(m-1)},\mj^{i-1},\mj \otimes_1 \beta,\vmj) \\
   &= (\mj^{(j-1)(m-1) +i-1},\iota \circ(\mj \otimes_1 \beta),\vmj) \overset{(2)}{=} (\mj^{(j-1)(m-1) +i-1},\mj,\vmj) = \varphi_{\Phi' (g,2)}.
\end{split}
\end{equation*}

\noindent 3.3. Suppose $\circ\Phi=\circ\Phi'$ is $d^n_j\circ
s^n_{j-1}=\mj$.

\vspace{1ex}

3.3.1. Suppose $g$ is $s^m_i$.
\begin{tabbing}
\hspace{3em}\=$\varphi_{\Phi (g,2)}$ \=$=(\mj^{(j-1)
m+i},\tau\circ(\kappa\otimes_1\mj),\vmj)$
\\[.5ex]
\>\>$=(\mj^{(j-1)m+i},\mj,\vmj)=\varphi_{\Phi' (g,2)}$,\quad by
(6).
\end{tabbing}

\vspace{1ex}

3.3.2. Suppose $g$ is $d^m_i$.
\begin{tabbing}
\hspace{3em}\=$\varphi_{\Phi (g,2)}$ \=$=(\mj^{(j-1)
(m-1)+i-1},\iota\circ(\beta\otimes_1\mj),\vmj)$
\\[.5ex]
\>\>$=(\mj^{(j-1)(m-1)+i-1},\mj,\vmj)=\varphi_{\Phi' (g,2)}$,\quad
by (3).
\end{tabbing}

\vspace{2ex}

\noindent 3.4. Suppose $\circ\Phi=\circ\Phi'$ is $d^n_j\circ
s^n_l=s^{n-1}_l\circ d^{n-1}_{j-1}$ for $j\geq l+2$.

\vspace{1ex}

3.4.1. Suppose $g$ is $s^m_i$.
\begin{tabbing}
\hspace{3em}\=$\varphi_{\Phi
(g,2)}=(\mj^{lm+i},\kappa,\mj^{(j-l-1)m-1},\tau,\vmj)=\varphi_{\Phi'
(g,2)}$.
\end{tabbing}

\vspace{1ex}

3.4.2. Suppose $g$ is $d^m_i$.
\begin{tabbing}
\hspace{3em}\=$\varphi_{\Phi
(g,2)}=(\mj^{l(m-1)+i-1},\beta,\mj^{(j-l-1)(m-1)-1},\iota,\vmj)=\varphi_{\Phi'
(g,2)}$.\`$\dashv$
\end{tabbing}

\prop{Lemma 4.4}{If $\Psi$ and $\Psi'$ are sequences of colour 1
such that $(\circ\Psi)^{\rm nf}$ is $\circ\Psi'$, and~$g$ is a
basic arrow of $\Delta^{\!op}$, then $\varphi_{\Psi
(g,2)}=\varphi_{\Psi' (g,2)}$.}

\dkz Let $\mu(\Psi)$ be a ``distance'' from $\circ\Psi$ to
$(\circ\Psi)^{\rm nf}$. For example, $\mu(\Psi)$ can be defined as
the ordered pair
\[
(n,m),
\]
where $n$ is the number of subsequences of $\Psi$ that are of the
form $(d,1)(s,1)$, i.e.,\ $s$ precedes $d$ looking from the right
to the left (not necessary immediately) in $\Psi$, and $m$ is the
number of subsequences of $\Psi$ of the form $(s_i,1)(s_j,1)$ with
$i\leq j$, or $(d_i,1)(d_j,1)$ with $i<j$. Suppose that our set of
``distances'' is lexicographically ordered.

We proceed by induction on $\mu(\Psi)$. If $\mu(\Psi)=(0,0)$, then
$\Psi=\Psi'$ and we are done. If $\mu(\Psi)>(0,0)$, then, by
Remark 3.2, $\Psi$ must be of the form $\Psi_2\Phi\Psi_1$, where
$\circ\Phi=\circ\Phi'$ is a basic equation of $\Delta^{\!op}$.
Then we have
\begin{tabbing}
\hspace{1.3em}\=$\varphi_{\Psi_2\Phi\Psi_1 (g,2)}$
\=$=\varphi_{\Psi_2 (g,2)}\WM_{\Phi\Psi_1}\,\circ\,
\WM_{\Psi_2}\,\varphi_{\Phi (g,2)}\WM_{\Psi_1}\,\circ\,
\WM_{\Psi_2\Phi}\,\varphi_{\Psi_1 (g,2)}$,
\\
\`(by Theorem 4.2)
\\[.5ex]
\>\>$=\varphi_{\Psi_2 (g,2)}\WM_{\Phi\Psi_1}\circ\,
\WM_{\Psi_2}\,\varphi_{\Phi' (g,2)}\WM_{\Psi_1}\circ\,
\WM_{\Psi_2\Phi'}\,\varphi_{\Psi_1 (g,2)}$,
\\
\`(by Lemma 4.3 and functoriality of $\WM_1$)
\\[.5ex]
\>\>$=\varphi_{\Psi_2\Phi'\Psi_1 (g,2)}$,\quad (by Theorem 4.2)
\\[.5ex]
\>\>$=\varphi_{\Psi' (g,2)}$.\quad (by the ind.\ hyp.\ since
$\mu(\Psi_2\Phi'\Psi_1)<\mu(\Psi_2\Phi\Psi_1)$)\` $\dashv$
\end{tabbing}

We can prove now $(i)$ by induction on the length of $\Gamma_1$
where in the induction step we use Lemma 4.4. We can prove $(ii)$
in a dual manner using the equations (7)-(12) for the proof of a
lemma dual to Lemma 4.3. So, we have:

\prop{Theorem 4.5}{The two-fold reduced bar construction $\WM$,
together with the natural transformations $\omega$, makes a lax
functor from $(\Delta^{\!op})^2$ to \emph{Cat}.}

\section{The three-fold monoidal categories}

The notion of three-fold monoidal category that we use in this
paper is defined in \cite[Section 7.1]{AM10} under the name
3-\textit{monoidal category}. In order to define this notion we
first define what the arrows between the two-fold monoidal
categories are.

\vspace{2ex}

\noindent{\sc Definition.}\quad A \emph{two-fold monoidal functor}
between two-fold monoidal categories $\cal C$ and $\cal D$ is a
5-tuple $\langle F,\sigma^1,\delta^1,\sigma^2,\delta^2\rangle$,
where for $i\in\{1,2\}$,
\[
\sigma^i_{A,B}\!:FA\otimes^{\cal D}_i FB\str F(A\otimes^{\cal C}_i
B)\quad{\rm and}\quad \delta^i\!:I^{\cal D}_i\str FI^{\cal C}_i
\]
are arrows of $\cal D$ natural in $A$ and $B$, such that $\langle
F,\sigma^1,\delta^1\rangle$ and $\langle
F,\sigma^2,\delta^2\rangle$ are monoidal functors between,
respectively, the first and the second monoidal structures of
$\cal C$ and $\cal D$. Moreover, the structure brought by the
arrows $\kappa$, $\beta$, $\tau$ and $\iota$ is preserved, which
means that the following four diagrams commute (with the
superscripts $\cal C$ and $\cal D$ omitted):

\begin{center}
\begin{picture}(80,50)
\put(0,40){\makebox(0,0){$I_1$}} \put(80,40){\makebox(0,0){$I_2$}}
\put(0,10){\makebox(0,0){$FI_1$}}
\put(80,10){\makebox(0,0){$FI_2$}}

\put(15,10){\vector(1,0){50}} \put(15,40){\vector(1,0){50}}
\put(0,33){\vector(0,-1){15}} \put(80,33){\vector(0,-1){15}}

\put(40,7){\makebox(0,0)[t]{$F\kappa$}}
\put(40,43){\makebox(0,0)[b]{$\kappa$}}
\put(-3,25){\makebox(0,0)[r]{$\delta^1$}}
\put(83,25){\makebox(0,0)[l]{$\delta^2$}}

\end{picture}
\end{center}

\begin{center}
\begin{picture}(280,80)
\put(0,70){\makebox(0,0){$I_1$}}
\put(80,70){\makebox(0,0){$I_1\otimes_2 I_1$}}
\put(0,10){\makebox(0,0){$FI_1$}}
\put(80,40){\makebox(0,0){$FI_1\otimes_2 FI_1$}}
\put(80,10){\makebox(0,0){$F(I_1\otimes_2 I_1)$}}

\put(15,10){\vector(1,0){35}} \put(15,70){\vector(1,0){40}}
\put(0,63){\vector(0,-1){45}} \put(80,63){\vector(0,-1){15}}
\put(80,33){\vector(0,-1){15}}

\put(35,7){\makebox(0,0)[t]{$F\beta$}}
\put(35,73){\makebox(0,0)[b]{$\beta$}}
\put(-3,40){\makebox(0,0)[r]{$\delta^1$}}
\put(83,55){\makebox(0,0)[l]{$\delta^1\otimes_2\delta^1$}}
\put(83,25){\makebox(0,0)[l]{$\sigma^2$}}

\put(280,70){\makebox(0,0){$I_2$}}
\put(200,70){\makebox(0,0){$I_2\otimes_1 I_2$}}
\put(280,10){\makebox(0,0){$FI_2$}}
\put(200,40){\makebox(0,0){$FI_2\otimes_1 FI_2$}}
\put(200,10){\makebox(0,0){$F(I_2\otimes_1 I_2)$}}

\put(225,10){\vector(1,0){35}} \put(225,70){\vector(1,0){40}}
\put(280,63){\vector(0,-1){45}} \put(200,63){\vector(0,-1){15}}
\put(200,33){\vector(0,-1){15}}

\put(245,7){\makebox(0,0)[t]{$F\tau$}}
\put(245,73){\makebox(0,0)[b]{$\tau$}}
\put(283,40){\makebox(0,0)[l]{$\delta^2$}}
\put(197,55){\makebox(0,0)[r]{$\delta^2\otimes_1\delta^2$}}
\put(197,25){\makebox(0,0)[r]{$\sigma^1$}}

\end{picture}
\end{center}

\begin{center}
\begin{picture}(200,80)

\put(0,70){\makebox(0,0){$(FA\otimes_2 FB)\otimes_1(FC\otimes_2
FD)$}} \put(200,70){\makebox(0,0){$(FA\otimes_1
FC)\otimes_2(FB\otimes_1 FD)$}}

\put(0,40){\makebox(0,0){$F(A\otimes_2 B)\otimes_1 F(C\otimes_2
D)$}} \put(200,40){\makebox(0,0){$F(A\otimes_1 C)\otimes_2
F(B\otimes_1 D)$}}

\put(0,10){\makebox(0,0){$F((A\otimes_2 B)\otimes_1(C\otimes_2
D))$}} \put(200,10){\makebox(0,0){$F((A\otimes_1
C)\otimes_2(B\otimes_1 D))$}}

\put(70,10){\vector(1,0){60}} \put(70,70){\vector(1,0){60}}
\put(0,63){\vector(0,-1){15}} \put(200,63){\vector(0,-1){15}}
\put(0,33){\vector(0,-1){15}} \put(200,33){\vector(0,-1){15}}

\put(100,7){\makebox(0,0)[t]{$F\iota$}}
\put(100,73){\makebox(0,0)[b]{$\iota$}}
\put(-3,55){\makebox(0,0)[r]{$\sigma^2\otimes_1\sigma^2$}}
\put(-3,25){\makebox(0,0)[r]{$\sigma^1$}}
\put(203,55){\makebox(0,0)[l]{$\sigma^1\otimes_2\sigma^1$}}
\put(203,25){\makebox(0,0)[l]{$\sigma^2$}}

\end{picture}
\end{center}

Let $Mon_2(Cat)$ be the 2-category whose 0-cells are the two-fold
monoidal categories, 1-cells are the two-fold monoidal functors,
and 2-cells are the \emph{two-fold monoidal transformations},
i.e.,\ monoidal transformations with respect to both the
structures. The monoidal structure of $Mon_2(Cat)$ is yet again
given by 2-products.

\vspace{2ex}

\noindent{\sc Definition.}\quad A \emph{three-fold monoidal}
category is a pseudomonoid in $Mon_2(Cat)$.

\vspace{2ex}

\noindent Hence, a three-fold monoidal category consists of the
following:

\vspace{.5ex}

1. a two-fold monoidal category $\cal M$,

\vspace{.5ex}

2. two-fold monoidal functors $\otm_3\!:{\cal M}\times {\cal
M}\str {\cal M}$ and $ I_3\!: 1\str {\cal M}$,

\vspace{.5ex}

3. two-fold monoidal transformations $\alpha_3$, $\rho_3$, and
$\lambda_3$ such that the structure $\langle{\cal M},\otm_3, I_3,
\alpha_3,\rho_3,\lambda_3\rangle$ satisfies the pseudomonoid
conditions.

\vspace{1ex}

In an unfolded form, this means that a three-fold monoidal
category is a category $\cal M$ equipped with three monoidal
structures ${\cal M}_1=\langle{\cal M},\otm_1, I_1\rangle$, ${\cal
M}_2=\langle{\cal M},\otm_2, I_2\rangle$, and ${\cal
M}_3=\langle{\cal M},\otm_3, I_3\rangle$ such that

\vspace{1ex}

\textbf{[1-2]} ${\cal M}_1$, ${\cal M}_2$, $\kappa\!: I_1\str
I_2$, $\beta\!: I_1\str I_1\otm_2 I_1$, $\tau\!: I_2\otm_1 I_2\str
I_2$, and $\iota\!:(A\otm_2 B)\otm_1(C\otm_2 D)\str(A\otm_1
C)\otm_2(B\otm_1 D)$,

\vspace{1ex}

\textbf{[2-3]} ${\cal M}_2$, ${\cal M}_3$, $\kappa'\!: I_2\str
I_3$, $\beta'\!: I_2\str I_2\otm_3 I_2$, $\tau'\!: I_3\otm_2
I_3\str I_3$, and $\iota'\!:(A\otm_3 B)\otm_2(C\otm_3
D)\str(A\otm_2 C)\otm_3(B\otm_2 D)$,

\vspace{1ex}

\textbf{[1-3]} ${\cal M}_1$, ${\cal M}_3$, $\kappa''\!: I_1\str
I_3$, $\beta''\!: I_1\str I_1\otm_3 I_1$, $\tau''\!: I_3\otm_1
I_3\str I_3$ and $\iota''\!:(A\otm_3 B)\otm_1(C\otm_3
D)\str(A\otm_1 C)\otm_3(B\otm_1 D)$,

\vspace{1ex}

\noindent are two-fold monoidal and, moreover, the following
equations hold:

\begin{tabbing}
\hspace{1.5em}\=(13)\hspace{1em}\=$\kappa'\circ\kappa=\kappa''$,
\\[1ex]
\>(14)\>$\beta'\circ\kappa=(\kappa\otimes_3\kappa)\circ\beta''$,
\\[1ex]
\>(15)\>$\tau'\circ(\kappa''\otimes_2\kappa'')\circ\beta=\kappa''$,
\\[1ex]
\>(16)\>$\iota'\circ(\beta''\otimes_2\beta'')\circ\beta=(\beta\otimes_3\beta)\circ\beta''$,
\\[1ex]
\>(17)\>$\kappa'\circ\tau=\tau''\circ(\kappa'\otimes_1\kappa')$,
\\[1ex]
\>(18)\>$\beta'\circ\tau=(\tau\otimes_3\tau)\circ\iota''\circ(\beta'\otimes_1\beta')$,
\\[1ex]
\>(19)\>$\tau'\circ(\tau''\otimes_2\tau'')\circ\iota=\tau''\circ(\tau'\otimes_1\tau')$,
\\[1ex]
\>(20)\>$\iota'\circ(\iota''\otimes_2\iota'')\circ\iota=(\iota\otimes_3\iota)\circ\iota''\circ(\iota'\otimes_1\iota')$.
\end{tabbing}
\textit{Note.} The last eight equations represent the four
commutative diagrams given above, with $F$ replaced by the
two-fold monoidal
functors $I_3$ and $\otimes_3$.\\

As in the case of two-fold monoidal categories, we are interested
only in \emph{three-fold strict monoidal} categories, i.e.,\ when
the structures ${\cal M}_1$, ${\cal M}_2$, and ${\cal M}_3$ are
strict monoidal.

The three-fold monoidal categories defined in
\cite[Definition~1.7]{BFSV} are the three-fold strict monoidal
categories from above in which, moreover, it is assumed that
$I_1=I_2=I_3=0$ and all the $\kappa$'s, $\beta$'s and $\tau$'s are
$\mj_0$. Hence, from the above list of eight equations, the
equations (13), (14), (15), and (17) are trivial, (16), (18) and
(19) are variants of internal and external unit conditions, while
the equation (20) corresponds to the \emph{big hexagonal
interchange diagram} (see \cite[Definition~1.7]{BFSV}).


\section{The three-fold reduced bar construction}

As in the two-fold case, we start with a definition of the
three-fold reduced bar construction based on a three-fold strict
monoidal category. Again, this construction corresponds to the one
given in the proof of \cite[Theorem~2.1]{BFSV}, save that the
latter construction is based on a category that is three-fold
monoidal in the sense of that paper.

For a three-fold strict monoidal category $\cal M$, we define
functions $\WM$ from objects and arrows of $(\Delta^{\!op})^3$ to
objects and arrows of \emph{Cat} in the following manner.

\vspace{2ex}

\noindent {\sc Definition.}\quad The \emph{three-fold reduced bar
construction} $\WM$ is defined on objects of $(\Delta^{\!op})^3$
as:
\[
\WM(n,m,p)=\cM^{n\cdot m\cdot p},
\]
and for arrows $f:n_s\str n_t$, $g:m_s\str m_t$ and $h:p_s\str
p_t$ of $\Delta^{\!op}$, we define $\WM(f,g,h)$ as the composition
\[
(\WM_3(h))^{n_t\cdot m_t}\circ(\WM_2^{p_s}(g))^{n_t}\circ
\WM_1^{m_s\cdot p_s}(f).
\]

For example, $\WM(d^2_1,s^2_1,d^2_1):{\cal M}^4\str {\cal M}^2$ is
defined as the composition
\[
(\WM_3(d^2_1))^2\circ \WM_2^{2}(s^2_1)\circ \WM_1^2(d^2_1),
\]
and for an object $(A,B,C,D)$ of ${\cal M}^4$ we have that
\[
\WM(d^2_1,s^2_1,d^2_1)(A,B,C,D)=((A\otimes_1
C)\otimes_3(B\otimes_1 D),I_2\otimes_3 I_2).
\]

As in the two-fold case, $\WM$ need not be a functor from
$(\Delta^{\!op})^3$ to \emph{Cat}, and our goal is to prove that
it is a lax functor. This means that for every composable pair of
arrows $e_1=(f_1,g_1,h_1)$ and  $e_2=(f_2,g_2,h_2)$ of
$(\Delta^{\!op})^3$, there is a natural transformation
\[
\omega_{e_2,e_1}\!:\WM(e_2)\circ \WM(e_1)\strt \WM(e_2\circ e_1),
\]
such that Diag~4.1 commutes.


We use coloured sequences and their shuffles in order to define
such natural transformations $\omega$. Let $\Phi$, $\Gamma$, and H
be sequences of colour 1, 2, and 3, respectively, such that
$\circ\Phi\!:n_s\str n_t$, $\circ\Gamma\!:m_s\str m_t$, and
$\circ{\rm H}\!:p_s\str p_t$. Let $\Theta$ be a shuffle of these
three sequences. For example, let $\Phi$ be $(d^2_2,1)(d^3_1,1)$,
let $\Gamma$ be $(d^2_1,2)$, let H be $(s^3_1,3)(s^2_1,3)$, and
let $\Theta$ be the following shuffle
\[
(d^2_2,1)(s^3_1,3)(d^3_1,1)(d^2_1,2)(s^2_1,3).
\]

For every member $(f,1)$ of $\Theta$, we define its \emph{inner
power} to be the product of the targets of its right-closest
$(g,2)$ and right-closest $(h,3)$ in $\Theta$. We may assume again
that such $(g,2)$ and $(h,3)$ exist since we can always add an
identity of colour 2 and an identity of colour 3 to the right of
$(f,1)$ in $\Theta$. For every member~$(g,2)$ of~$\Theta$, we
define its \emph{inner power} to be the target of its
right-closest $(h,3)$ in $\Theta$, and we define its \emph{outer
power} to be the target of its right-closest $(f,1)$ in~$\Theta$.
For every member $(h,3)$ of $\Theta$, we define its \emph{outer
power} to be the product of the targets of its right-closest
$(f,1)$ and right-closest $(g,2)$ in $\Theta$. For $\Theta$ as
above, for example, we have that the outer power of $(s^2_1,3)$
is~6.

Let $\cal M$ be a three-fold strict monoidal category. We define a
functor
\[
\WM_\Theta\!:{\cal M}^{n_s\cdot m_s\cdot p_s}\str{\cal
M}^{n_t\cdot m_t\cdot p_t}
\]
in the following way: replace in $\Theta$ each $(f,1)$ whose inner
power is $i$ by $\WM_1^i(f)$, every $(g,2)$ whose inner power is
$i$ and outer power is $o$ by $(\WM_2^i(g))^o$ and every $(h,3)$
whose outer power is $o$ by $(\WM_3(h))^o$, and insert $\circ$'s.
For $\Theta$ as above, we have that $\WM_\Theta$ is
\[
\WM_1^3(d^2_2)\circ (\WM_3(s^3_1))^2 \circ \WM_1^2(d^3_1)\circ
(\WM_2^2(d^2_1))^3 \circ (\WM_3(s^2_1))^6,
\]
which gives that $\WM_\Theta(A,B,C,D,E,F)$ is the 3-tuple
\[
((A\otimes_2 B)\otimes_1(C\otimes_2 D),I_3,(I_3\otimes_2
I_3)\otimes_1(I_3\otimes_2 I_3)).\]


It is easy to see that for arrows $f$, $g$ and $h$ of
$\Delta^{\!op}$, we have that
\[
\WM(f,g,h)=\WM_{{\rm
H}\Gamma\Phi},
\]
for arbitrary $\Phi$ of colour~1, $\Gamma$ of colour~2 and H of
colour 3, such that $\circ\Phi=f$, $\circ\Gamma=g$ and $\circ{\rm
H}=h$. This may serve as a \emph{combinatorial} definition of the
three-fold reduced bar construction $\WM$ (cf.\ the combinatorial
definition of the two-fold reduced bar construction given in
Section~4).


For basic arrows $f\!:n\str n'$, $g\!:m\str m'$ of $\Delta^{\!op}$
and $w\geq 0$ we define a natural transformation
\[
\chi^{1,2}_w(f,g)\!:\WM_1^{m'\!w}(f)\circ(\WM_2^w(g))^n\;\strt\;
(\WM_2^w(g))^{n'}\circ\WM_1^{mw}(f)
\]
to be the identity natural transformation except in the following
cases:

\begin{center}
\begin{tabular}{c|c|c}
$f$ & $g$ & $\chi^{1,2}_w(f,g)$
\\
\hline\hline\\[-2.2ex] $s^{n+1}_j$ & $s^{m+1}_i$ &
$\displaystyle(\mj^{j(m+1)w},\underbrace{\mj^{iw},\kappa^w,\vmj}_{(m+1)w},\vmj)$
\\
\hline\\[-2.2ex] $d^n_j$, $1\leq j\leq n-1$ & $s^{m+1}_i$ &
$\displaystyle(\mj^{(j-1)(m+1)w},\underbrace{\mj^{iw},\tau^w,\vmj}_{(m+1)w},\vmj)$
\\
\hline\\[-2.2ex] $s^{n+1}_j$ & $d^m_i$, $1\leq i\leq m-1$ &
$\displaystyle(\mj^{j(m-1)w},\underbrace{\mj^{(i-1)w},\beta^w,\vmj}_{(m-1)w},\vmj)$
\\
\hline\\[-2.2ex] $d^n_j$, $1\leq j\leq n-1$ & $d^m_i$, $1\leq i\leq m-1$ &
$\displaystyle(\mj^{(j-1)(m-1)w},\underbrace{\mj^{(i-1)w},\iota^w,\vmj}_{(m-1)w},\vmj)$
\end{tabular}
\captionof{table}{$\chi_w^{1,2}$ in nontrivial cases.}
\label{table:3dim12}
\end{center}

\vspace{1ex}

In order to simplify some calculations and improve the
presentation of the paper, we introduce the following formal
operation of \emph{multiplication} (always from the right) of
tuples representing the natural transformations by 0-1 matrices
having in each column exactly one entry equal to 1 and all the
other entries equal to 0, which is derived from the standard
multiplication of matrices. For example,
\[
(\mj,\kappa,\mj)\le
\begin{array}{ccccccccccc} 1  & 0 & 0 & 0 & 0 & 1 & 1 & 0 & 0 & 0 & 0 \\
0 & 1 & 1 & 0 & 0 & 0 & 0 & 1 & 1 & 0 & 0 \\ 0 & 0 & 0 & 1 & 1 & 0
& 0 & 0 & 0 & 1 & 1
\end{array} \re =(\mj,\kappa^2,\mj^2, \mj^2,\kappa^2,\mj^2).
\]
Note that the tuples of the third column of Table
\ref{table:3dim12} are obtained as a result of multiplication of
the tuples in the third column of Table~\ref{table:2dim} by the
matrix
\[I_{n'}\otimes I_{m'}\otimes(\underbrace{1,\ldots,1}_w),\] where
$I_k$ is the $k\times k$  identity matrix and $\otimes$ is the
\emph{Kronecker product} of matrices.

For basic arrows $g\!:m\str m'$, $h\!:p\str p'$ of $\Delta^{\!op}$
and $u\geq 0$ we define a natural transformation
\[
\chi^{2,3}_u(g,h)\!:(\WM_2^{p'}(g))^u\circ(\WM_3(h))^{um}\;\strt\;
(\WM_3(h))^{um'}\circ(\WM_2^p(g))^u
\]
to be the identity natural transformation except in the following
cases: \vspace*{1ex}
\begin{center}
\begin{tabular}{c|c|c}
$g$ & $h$ & $\chi^{2,3}_u(g,h)$
\\
\hline\hline\\[-2.2ex] $s^{m+1}_i$ & $s^{p+1}_k$ &
$\displaystyle((\underbrace{\mj^{i(p+1)},\mj^k,\kappa',\vmj}_{(m+1)(p+1)})^u)$
\\
\hline\\[-2.2ex] $d^m_i$, $1\leq i\leq m-1$ & $s^{p+1}_k$ &
$\displaystyle((\underbrace{\mj^{(i-1)(p+1)},\mj^k,\tau',\vmj}_{(m-1)(p+1)})^u)$
\\
\hline\\[-2.2ex] $s^{m+1}_i$ & $d^p_k$, $1\leq k\leq p-1$ &
$\displaystyle((\underbrace{\mj^{i(p-1)},\mj^{k-1},\beta',\vmj}_{(m+1)(p-1)})^u)$
\\
\hline\\[-2.2ex] $d^m_i$, $1\leq i\leq m-1$ & $d^p_k$, $1\leq k\leq p-1$ &
$\displaystyle((\underbrace{\mj^{(i-1)(p-1)},\mj^{k-1},\iota',\vmj}_{(m-1)(p-1)})^u)$
\end{tabular}
\captionof{table}{$\chi_u^{2,3}$ in nontrivial cases.}
\label{table:3dim23}
\end{center}

\vspace{1ex}

\noindent Note that the tuples of the third column of this table
are obtained as a result of multiplication of the tuples in the
third column of Table~\ref{table:2dim} (where $m$ is replaced by
$p$, $n$ is replaced by $m$, $i$ is replaced by $k$, $j$ is
replaced by $i$, and $\kappa$, $\beta$, $\tau$, and $\iota$ are
replaced by $\kappa'$, $\beta'$, $\tau'$, and $\iota'$) by the
matrix
\[(\underbrace{1,\ldots,1}_u)\otimes I_{m'}\otimes
I_{p'}.\]

Finally, for basic arrows $f\!:n\str n'$, $h\!:p\str p'$ of
$\Delta^{\!op}$ and $v\geq 0$ we define a natural transformation
\[
\chi^{1,3}_v(f,h)\!:\WM_1^{vp'}(f)\circ(\WM_3(h))^{nv}\;\strt\;
(\WM_3(h))^{n'\!v}\circ\WM_1^{vp}(f)
\]
to be the identity natural transformation except in the following
cases: \vspace{1ex}
\begin{center}
\begin{tabular}{c|c|c}
$f$ & $h$ & $\chi^{1,3}_v(f,h)$
\\
\hline\hline\\[-2.2ex] $s^{n+1}_j$ & $s^{p+1}_k$ &
$\displaystyle(\mj^{jv(p+1)},(\underbrace{\mj^{k},\kappa'',\vmj}_{p+1})^v,\vmj)$
\\
\hline\\[-2.2ex] $d^n_j$, $1\leq j\leq n-1$ & $s^{p+1}_k$ &
$\displaystyle(\mj^{(j-1)v(p+1)},(\underbrace{\mj^k,\tau'',\vmj}_{p+1})^v,\vmj)$
\\
\hline\\[-2.2ex] $s^{n+1}_j$ & $d^p_k$, $1\leq k\leq p-1$ &
$\displaystyle(\mj^{jv(p-1)},(\underbrace{\mj^{k-1},\beta'',\vmj}_{p-1})^v,\vmj)$
\\
\hline\\[-2.2ex] $d^n_j$, $1\leq j\leq n-1$ & $d^p_k$, $1\leq k\leq p-1$ &
$\displaystyle(\mj^{(j-1)v(p-1)},(\underbrace{\mj^{k-1},\iota'',\vmj}_{p-1})^v,\vmj)$
\end{tabular}
\captionof{table}{$\chi_v^{1,3}$ in nontrivial cases.}
\label{table:3dim13}
\end{center}

\vspace{1ex}

\noindent As in the previous cases, the tuples of the third column
of this table are obtained as a result of multiplication of the
tuples in the third column of Table~\ref{table:2dim} (with some
necessary replacements) by a certain matrix, in this case that
matrix is
\[I_{n'}\otimes(\underbrace{1,\ldots,1}_v)\otimes
I_{p'}.\]

\prop{Lemma 6.1}{For basic arrows $f\!:n\str n'$, $g\!:m\str m'$,
and $h\!:p\str p'$ of $\Delta^{\!op}$ the following diagram
commutes:}

\begin{center}
\begin{picture}(200,135)

\put(100,130){\makebox(0,0){$\WM_{(f,1)(g,2)(h,3)}$}}
\put(30,90){\makebox(0,0){$\WM_{(g,2)(f,1)(h,3)}$}}
\put(170,90){\makebox(0,0){$\WM_{(f,1)(h,3)(g,2)}$}}
\put(30,50){\makebox(0,0){$\WM_{(g,2)(h,3)(f,1)}$}}
\put(170,50){\makebox(0,0){$\WM_{(h,3)(f,1)(g,2)}$}}
\put(100,10){\makebox(0,0){$\WM_{(h,3)(g,2)(f,1)}$}}

\put(70,120){\vector(-3,-2){30}} \put(130,120){\vector(3,-2){30}}
\put(30,80){\vector(0,-1){23}} \put(170,80){\vector(0,-1){23}}
\put(40,40){\vector(3,-2){30}} \put(160,40){\vector(-3,-2){30}}

\put(45,115){\makebox(0,0)[r]{$\chi^{1,2}_{p'}(f,g)(\WM_3(h))^{nm}$}}

\put(25,70){\makebox(0,0)[r]{
$(\WM_2^{p'}(g))^{n'}\chi^{1,3}_m(f,h)$}}

\put(45,25){\makebox(0,0)[r]{$\chi^{2,3}_{n'}(g,h)\WM_1^{pm}(f)$}}

\put(155,115){\makebox(0,0)[l]{
$\WM_1^{p'm'}(f)\chi^{2,3}_n(g,h)$}}

\put(175,70){\makebox(0,0)[l]{$\chi^{1,3}_{m'}(f,h)(\WM_2^p(g))^n$}}

\put(155,25){\makebox(0,0)[l]{
$(\WM_3(h))^{n'm'}\chi^{1,2}_p(f,g)$}}

\end{picture}
\end{center}

\dkz Consider the following table in which $d^x_y$ is such that
$0<y<x$.

\begin{center}
\begin{tabular}{l|l|l|c|c}
$f$ & $g$ & $h$ & component & EQ
\\
\hline\hline\\[-2.2ex] $s^{n+1}_j$ & $s^{m+1}_i$ & $s^{p+1}_k$ &
$j(m+1)(p+1)+i(p+1)+k+1$ & (13)
\\
\hline\\[-2.2ex] $s^{n+1}_j$ & $s^{m+1}_i$ & $d^p_k$ &
$j(m+1)(p-1)+i(p-1)+k$ & (14)
\\
\hline\\[-2.2ex] $s^{n+1}_j$ & $d^m_i$ & $s^{p+1}_k$ &
$j(m-1)(p+1)+(i-1)(p+1)+k+1$ & (15)
\\
\hline\\[-2.2ex] $s^{n+1}_j$ & $d^m_i$ & $d^p_k$ &
$j(m-1)(p-1)+(i-1)(p-1)+k$ & (16)
\\
\hline\\[-2.2ex] $d^n_j$ & $s^{m+1}_i$ & $s^{p+1}_k$ &
$(j-1)(m+1)(p+1)+i(p+1)+k+1$ & (17)
\\
\hline\\[-2.2ex] $d^n_j$ & $s^{m+1}_i$ & $d^p_k$ &
$(j-1)(m+1)(p-1)+i(p-1)+k$ & (18)
\\
\hline\\[-2.2ex] $d^n_j$ & $d^m_i$ & $s^{p+1}_k$ &
$(j-1)(m-1)(p+1)+(i-1)(p+1)+k+1$ & (19)
\\
\hline\\[-2.2ex] $d^n_j$ & $d^m_i$ & $d^p_k$ &
$(j-1)(m-1)(p-1)+(i-1)(p-1)+k$ & (20)
\end{tabular}
\end{center}

This gives a list of all nontrivial cases for $f$, $g$, and $h$.
In this table we point out the component of the two $n'\cdot
m'\cdot p'$-tuples of arrows, representing the left-hand side and
the right-hand side of the above diagram, where we use one of the
equations (13)-(20). In all the other components, the left-hand
side is equal to the right-hand side by simple categorial
arguments.

As an illustration of these arguments, here we give a proof for
one of the cases from the table, namely when $f=d_j^n$,
$g=s_i^{m+1}$, and $h=d^p_k$. At the left hand side of the diagram
we have the following
\begin{align*}
&\chi^{1,2}_{p-1}(d_j^n,s_i^{m+1})(\WM_3(d^p_k))^{nm} =
(\mj^{(j-1)(m+1)(p-1)}, \underbrace{\mj^{i(p-1)},
\tau^{p-1},\vmj}_{\scriptscriptstyle(m+1)(p-1)},\vmj),
\label{diag:3dimL1}
\tag{L1} \\
 &(\WM_2^{p-1}(s_i^{m+1}))^{n-1}\chi^{1,3}_m(d_j^n,d^p_k) = \\&\hphantom{=======}=(\WM_2^{p-1}(s_i^{m+1}))^{n-1}\; (\mj^{(j-1)m(p-1)},(\underbrace{\mj^{k-1},\iota'',\vmj}_{\scriptscriptstyle p-1})^m,\vmj) \\
 &\hphantom{=======}= (\mj^{(j-1)(m+1)(p-1)},(\underbrace{\mj^{k-1},\iota'',\vmj}_{\scriptscriptstyle p-1})^i,\mj^{\scriptscriptstyle p-1},(\underbrace{\mj^{k-1},\iota'',\vmj}_{\scriptscriptstyle p-1})^{m-i},\vmj),
 \tag{L2}\label{diag:3dimL2} \\
 &\chi^{2,3}_{n-1}(s_i^{m+1},d^p_k)\WM_1^{mp}(d_j^n) = ((\underbrace{\mj^{i(p-1)},\mj^{k-1},\beta',\vmj}_{\scriptscriptstyle (m+1)(p-1)})^{n-1}), \tag{L3}\label{diag:3dimL3}
\end{align*}
while at the right hand side we have:
\begin{align*}
 &\WM_1 ^{(m+1)(p-1)} (d_j^n)\chi^{2,3}_n(s_i^{m+1},d^p_k) = \WM_1^{(m+1)(p-1)}(d_j^n) ((\underbrace{\mj^{i(p-1)},\mj^{k-1},\beta',\vmj}_{\scriptscriptstyle (m+1)(p-1)})^n) \\
 &=  ((\underbrace{\mj^{i(p-1)},\mj^{k-1},\beta',\vmj}_{\scriptscriptstyle (m+1)(p-1)})^{j-1}, \underbrace{\mj^{i(p-1)+k-1},\beta'\otimes_1 \beta',\vmj}_{\scriptscriptstyle (m+1)(p-1)},(\underbrace{\mj^{i(p-1)},\mj^{k-1},\beta',\vmj}_{\scriptscriptstyle (m+1)(p-1)})^{n-j-1}),
 \tag{D1}
 \label{diag:3dimD1}\\
 &\chi^{1,3}_{m+1}(d_j^n,d^p_k)(\WM_2^p(s_i^{m+1}))^n = (\mj^{(j-1)(m+1)(p-1)},(\underbrace{\mj^{k-1},\iota'',\vmj}_{\scriptscriptstyle p-1})^{m+1},\vmj),
 \tag{D2}
 \label{diag:3dimD2}\\
 &(\WM_3(d^p_k))^{(n-1)(m-1)}\chi^{1,2}_p(d_j^n,s_i^{m+1}) = \\
 &\hphantom{=======}=(\WM_3(d^p_k))^{(n-1)(m-1)}\circ (\mj^{(j-1)(m+1)p},\underbrace{\mj^{ip},\tau^p,\vmj}_{\scriptscriptstyle (m+1)p},\vmj)\\
 &\hphantom{=======}=(\mj^{(j-1)(m+1)(p-1)},\underbrace{\mj^{i(p-1)},\tau^{k-1},\tau \otimes_3 \tau, \tau^{p-1-k},\vmj}_{\scriptscriptstyle (m+1)(p-1)},\vmj).
 \tag{D3}
 \label{diag:3dimD3}
\end{align*}
Now we take a look at all entries that are not equal to $\mj$
(non-identities). For example, in (\ref{diag:3dimL1}) the
non-identities are at positions
\[
(j-1)(m+1)(p+1)+i(p-1)+l,\quad {\rm for}\; 1\leq l \leq p-1,
\]
and those entries are equal to $\tau$. By comparing the
non-identities for (\ref{diag:3dimD1}), (\ref{diag:3dimD2}),
(\ref{diag:3dimD3}), (\ref{diag:3dimL1}), (\ref{diag:3dimL2}), and
(\ref{diag:3dimL3}), we get that the only difference is at
position $(j-1)(m+1)(p+1)+i(p-1)+k$, where we have that
$\beta'\circ\mj\circ \tau$ must be equal to $(\tau\otimes_3
\tau)\circ \iota''\circ(\beta'\otimes_1 \beta')$, which is exactly
our equation~(18). \qed

Let $\Theta_0,\ldots,\Theta_j$ for $j\geq 0$ be shuffles of
$\Phi$, $\Gamma$, and H such that $\Theta_0=\Theta$ and
$\Theta_j={\rm H}\Gamma\Phi$, and if $j>0$, then for every $0\leq
i\leq j-1$ we have that $\Theta_i=\Pi(f,1)(g,2)\Lambda$ and
$\Theta_{i+1}=\Pi(g,2)(f,1)\Lambda$, or
$\Theta_i=\Pi(g,2)(h,3)\Lambda$ and
$\Theta_{i+1}=\Pi(h,3)(g,2)\Lambda$, or
$\Theta_i=\Pi(f,1)(h,3)\Lambda$ and
$\Theta_{i+1}=\Pi(h,3)(f,1)\Lambda$. We call
$\Theta_0,\ldots,\Theta_j$ a \emph{normalizing path} starting with
$\Theta$. Its \emph{length} is $j$ and Proposition 4.1 still
holds.

If $\Theta_i=\Pi(f,1)(g,2)\Lambda$ and
$\Theta_{i+1}=\Pi(g,2)(f,1)\Lambda$, then for $w$ being the target
of the leftmost member of $\Lambda$ of colour 3 we have that
\[
\varphi_i=\WM_\Pi\;\chi^{1,2}_w(f,g)\;\WM_\Lambda
\]
is a natural transformation from $\WM_{\Theta_i}$ to
$\WM_{\Theta_{i+1}}$. We define $\varphi_i$ analogously in the
other two possibilities for the pair $\Theta_i$, $\Theta_{i+1}$
relying on $\chi^{2,3}_u(g,h)$ or $\chi^{1,3}_v(f,h)$, for $u$
being the target of the leftmost member of $\Lambda$ of colour 1
and $v$ being the target of the leftmost member of $\Lambda$ of
colour 2. We define $\varphi(\Theta_0,\ldots,\Theta_k)$ as in the
two-fold case and for $\Theta'_0,\ldots,\Theta'_k$ being another
normalizing path starting with $\Theta$, we can show the
following.

\prop{Theorem 6.2}{
$\varphi(\Theta_0,\ldots,\Theta_j)=\varphi(\Theta'_0,\ldots,\Theta'_j)$.}

\dkz We proceed by induction on $j\geq 0$. If $j=0$, then
$\varphi(\Theta_0)=\varphi(\Theta'_0)=\mj$.

If $j>0$, then we are either in the situation as in the proof of
Theorem 4.2 and we proceed analogously, or for some basic arrows
$f\!:n\str n'$, $g\!:m\str m'$, and $h\!:p\str p'$ of
$\Delta^{\!op}$ we have that
\[
\varphi_0=\WM_\Pi\;\chi^{1,2}_{p'}(f,g)\;\WM_{(h,3)\Lambda}\quad
{\rm and}\quad
\varphi'_0=\WM_{\Pi(f,1)}\;\chi^{2,3}_n(g,h)\;\WM_{\Lambda}.
\]
In the latter case, we use Lemma 6.1 and the induction hypothesis
twice to obtain the following commutative diagram.

\begin{center}
\begin{picture}(200,180)(0,-40)

\put(100,130){\makebox(0,0){$\WM_\Pi\circ\WM_{(f,1)(g,2)(h,3)}\circ\WM_\Lambda$}}
\put(15,90){\makebox(0,0){$\WM_\Pi\circ\WM_{(g,2)(f,1)(h,3)}\circ\WM_\Lambda$}}
\put(185,90){\makebox(0,0){$\WM_\Pi\circ\WM_{(f,1)(h,3)(g,2)}\circ\WM_\Lambda$}}
\put(15,50){\makebox(0,0){$\WM_\Pi\circ\WM_{(g,2)(h,3)(f,1)}\circ\WM_\Lambda$}}
\put(185,50){\makebox(0,0){$\WM_\Pi\circ\WM_{(h,3)(f,1)(g,2)}\circ\WM_\Lambda$}}
\put(100,10){\makebox(0,0){$\WM_\Pi\circ\WM_{(h,3)(g,2)(f,1)}\circ\WM_\Lambda$}}
\put(100,-35){\makebox(0,0){$\WM_{{\rm H}\Gamma\Phi}$}}

\put(70,120){\vector(-3,-2){30}} \put(130,120){\vector(3,-2){30}}
\put(30,80){\vector(0,-1){23}} \put(170,80){\vector(0,-1){23}}
\put(40,40){\vector(3,-2){30}} \put(160,40){\vector(-3,-2){30}}
\put(100,-10){\vector(0,-1){15}} \put(30,-35){\vector(1,0){45}}
\put(170,-35){\vector(-1,0){45}}

\put(-40,80){\vector(-3,-2){30}} \put(240,80){\vector(3,-2){30}}

\put(45,115){\makebox(0,0)[r]{$\varphi_0$}}
\put(-60,75){\makebox(0,0)[r]{ $\varphi_1$}}
\put(50,-30){\makebox(0,0)[b]{ $\varphi_{j-1}$}}

\put(155,115){\makebox(0,0)[l]{ $\varphi_0'$}}
\put(260,75){\makebox(0,0)[l]{ $\varphi_1'$}}
\put(150,-30){\makebox(0,0)[b]{ $\varphi_{j-1}'$}}

\put(-20,20){\makebox(0,0){\scriptsize ind. hyp.}}
\put(220,20){\makebox(0,0){\scriptsize ind. hyp.}}

\put(100,70){\makebox(0,0){\scriptsize Lemma 6.1}}

\put(100,0){\makebox(0,0){$\vdots$}}

\qbezier[60](-75,55)(-100,-10)(20,-35)
\qbezier[60](275,55)(300,-10)(180,-35)

\put(270,-35){\makebox(0,0){$\dashv$}}

\end{picture}
\end{center}

\vspace{1ex}

By Theorem 6.2, the following definition is correct.

\vspace{2ex}

\noindent {\sc Definition.}\quad Let $\varphi_\Theta:\WM_\Theta
\strt \WM_{{\rm H}\Gamma\Phi}$ be
$\varphi(\Theta_0,\ldots,\Theta_j)$, for an arbitrary normalizing
path $\Theta_0,\ldots,\Theta_j$ starting with $\Theta$.


\vspace{2ex}

We are ready to define a natural transformation
\[
\omega_{e_2,e_1}\!:\WM(e_2)\circ \WM(e_1)\strt \WM(e_2\circ e_1),
\]
for every composable pair of arrows $e_1=(f_1,g_1,h_1)$ and
$e_2=(f_2,g_2,h_2)$ of $(\Delta^{\!op})^3$.

\vspace{2ex}

\noindent {\sc Definition.}\quad Let $\Phi_1$ and $\Phi_2$ be
sequences of colour 1, let $\Gamma_1$ and $\Gamma_2$ be sequences
of colour 2, and let ${\rm H}_1$ and ${\rm H}_2$ be sequences of
colour 3, such that $\circ\Phi_1$ is $f^{\rm nf}_1$, $\circ\Phi_2$
is $f^{\rm nf}_2$, $\circ\Gamma_1$ is $g^{\rm nf}_1$,
$\circ\Gamma_2$ is $g^{\rm nf}_2$, $\circ{\rm H}_1$ is $h^{\rm
nf}_1$, and $\circ{\rm H}_2$ is $h^{\rm nf}_2$. We define
\[
\omega_{e_2,e_1}\quad {\rm as}\quad \varphi_{{\rm
H}_2\Gamma_2\Phi_2{\rm H}_1\Gamma_1\Phi_1}.
\]

It remains to prove that our Diag~4.1 commutes. Let
$e_1=(f_1,g_1,h_1)$, $e_2=(f_2,g_2,h_2)$ and $e_3=(f_3,g_3,h_3)$
be such that the composition $e_3\circ e_2\circ e_1$ is defined in
$(\Delta^{\!op})^3$. Let $\Phi_1$, $\Phi_2$, $\Gamma_1$,
$\Gamma_2$, ${\rm H}_1$ and ${\rm H}_2$ be as above, and let
$\Phi_3$, $\Gamma_3$ and ${\rm H}_3$ be sequences of colour~1, 2
and 3 respectively such that $\circ\Phi_3$ is~$f^{\rm nf}_3$,
$\circ\Gamma_3$ is~$g^{\rm nf}_3$ and $\circ{\rm H}_3$ is~$h^{\rm
nf}_3$. In this case, Diag~4.1 reads

\begin{center}
\begin{picture}(100,75)(0,5)

\put(50,75){\makebox(0,0){$\WM(e_3)\circ\WM(e_2)\circ \WM(e_1)$}}
\put(-40,40){\makebox(0,0){$\WM(e_3\circ e_2)\circ \WM(e_1)$}}
\put(140,40){\makebox(0,0){$\WM(e_3)\circ\WM(e_2\circ e_1)$}}
\put(50,10){\makebox(0,0){$\WM(e_3\circ e_2\circ e_1)$}}

\put(-10,60){\vector(-2,-1){20}} \put(110,60){\vector(2,-1){20}}
\put(-30,30){\vector(2,-1){30}} \put(130,30){\vector(-2,-1){30}}

\put(-20,60){\makebox(0,0)[r]{\scriptsize $\varphi_{{\rm
H}_3\Gamma_3\Phi_3 {\rm H}_2\Gamma_2\Phi_2}\WM_{{\rm
H}_1\Gamma_1\Phi_1}$}} \put(120,60){\makebox(0,0)[l]{\scriptsize
$\WM_{{\rm H}_3\Gamma_3\Phi_3} \varphi_{{\rm
H}_2\Gamma_2\Phi_2{\rm H}_1\Gamma_1\Phi_1}$}}

\put(-20,20){\makebox(0,0)[r]{\scriptsize $\varphi_{{\rm
H}'\Gamma'\Phi'{\rm H}_1\Gamma_1\Phi_1}$}}
\put(120,20){\makebox(0,0)[l]{\scriptsize $\varphi_{{\rm
H}_3\Gamma_3\Phi_3{\rm H}''\Gamma''\Phi''}$}}

\end{picture}
\end{center}
where $\circ\Phi'$ is $(f_3\circ f_2)^{\rm nf}$, $\circ\Gamma'$ is
$(g_3\circ g_2)^{\rm nf}$, $\circ{\rm H}'$ is $(h_3\circ h_2)^{\rm
nf}$, $\circ\Phi''$ is $(f_2\circ f_1)^{\rm nf}$, $\circ\Gamma''$
is $(g_2\circ g_1)^{\rm nf}$ and $\circ{\rm H}''$ is $(h_2\circ
h_1)^{\rm nf}$.

\vspace{2ex}

By Theorem 6.2 we have the following commutative diagram

\begin{center}
\begin{picture}(100,75)(0,5)

\put(50,75){\makebox(0,0){$\WM(e_3)\circ\WM(e_2)\circ \WM(e_1)$}}
\put(-40,40){\makebox(0,0){$\WM(e_3\circ e_2)\circ \WM(e_1)$}}
\put(140,40){\makebox(0,0){$\WM(e_3)\circ\WM(e_2\circ e_1)$}}
\put(50,10){\makebox(0,0){$\WM(e_3\circ e_2\circ e_1)$}}

\put(-10,60){\vector(-2,-1){20}} \put(110,60){\vector(2,-1){20}}
\put(-30,30){\vector(2,-1){30}} \put(130,30){\vector(-2,-1){30}}
\put(50,65){\line(0,-1){5}} \put(50,50){\vector(0,-1){30}}

\put(-20,60){\makebox(0,0)[r]{\scriptsize $\varphi_{{\rm
H}_3\Gamma_3\Phi_3 {\rm H}_2\Gamma_2\Phi_2}\WM_{{\rm
H}_1\Gamma_1\Phi_1}$}} \put(120,60){\makebox(0,0)[l]{\scriptsize
$\WM_{{\rm H}_3\Gamma_3\Phi_3} \varphi_{{\rm
H}_2\Gamma_2\Phi_2{\rm H}_1\Gamma_1\Phi_1}$}}

\put(-20,20){\makebox(0,0)[r]{\scriptsize $\varphi_{{\rm H}_3{\rm
H}_2\Gamma_3\Gamma_2\Phi_3\Phi_2{\rm H}_1\Gamma_1\Phi_1}$}}
\put(120,20){\makebox(0,0)[l]{\scriptsize $\varphi_{{\rm
H}_3\Gamma_3\Phi_3{\rm H}_2{\rm
H}_1\Gamma_2\Gamma_1\Phi_2\Phi_1}$}}
\put(50,55){\makebox(0,0){\scriptsize $\varphi_{{\rm
H}_3\Gamma_3\Phi_3{\rm H}_2\Gamma_2\Phi_2{\rm
H}_1\Gamma_1\Phi_1}$}}

\end{picture}
\end{center}
Hence, to prove that Diag~4.1 commutes, it suffices to show that


\begin{tabbing}
\hspace{1.5em}\=$(i)$\quad\=$\varphi_{{\rm H}_3{\rm
H}_2\Gamma_3\Gamma_2\Phi_3\Phi_2{\rm
H}_1\Gamma_1\Phi_1}=\varphi_{{\rm H}'\Gamma'\Phi'{\rm
H}_1\Gamma_1\Phi_1}$\quad  and
\\[2ex]
\>$(ii)$\>$\varphi_{{\rm H}_3\Gamma_3\Phi_3{\rm H}_2{\rm
H}_1\Gamma_2\Gamma_1\Phi_2\Phi_1}=\varphi_{{\rm
H}_3\Gamma_3\Phi_3{\rm H}''\Gamma''\Phi''}$.
\end{tabbing}

To prove $(i)$ and $(ii)$ we use the same arguments as in the
two-fold case. Let $x$, $y$, and $z$ be three different elements
of the set $\{1,2,3\}$ such that $x<y$. Note that the position of
$(\mj_q,z)$ in the two shuffles of the lemma below is irrelevant;
$(\mj_q,z)$ serves just to keep $\varphi$ correctly defined and to
introduce the parameter~$q$.

\prop{Lemma 6.3}{If $\Phi$ and $\Phi'$ are sequences of colour $x$
such that $\circ\Phi=\circ\Phi'$ is a basic equation of
$\Delta^{\!op}$, and $g$ is a basic arrow of $\Delta^{\!op}$, then
for every $q\geq 0$ we have that $\varphi_{\Phi (g,y)
(\mjs_q,z)}=\varphi_{\Phi' (g,y)(\mjs_q,z)}$.}

\dkz Suppose the target of $\circ\Phi$ is $n'$ and the target of
$g$ is $m'$. If $x=1$, $y=2$, and $z=3$, then we proceed as in
Lemma 4.3 with all the cases modified so that the tuples
representing the natural transformations are multiplied by the
matrix $\displaystyle I_{n'}\otimes
I_{m'}\otimes(\underbrace{1,\ldots,1}_q)$. For example, Case 1.1.1
now reads
\begin{tabbing}
\hspace{1.5em}\=$\varphi_{\Phi (g,2)(\mjs_q,3)}=(\mj^{(j-1)
mq},\mj^{iq},\tau^q,\mj^{((l-j-1)m-1)q},\tau^q,\vmj)=\varphi_{\Phi'
(g,2)(\mjs_q,3)}$.
\end{tabbing}

If $x=2$, $y=3$, and $z=1$, we again proceed as in Lemma 4.3 with
all the cases modified so that $\kappa$, $\beta$, $\tau$, and
$\iota$ are replaced by $\kappa'$, $\beta'$, $\tau'$, and
$\iota'$, and the tuples representing the natural transformations
are multiplied by the matrix $\displaystyle
(\underbrace{1,\ldots,1}_q) \otimes I_{n'}\otimes I_{m'}$. For
example, Case 1.1.1 now reads
\begin{tabbing}
\hspace{1.5em}\=$\varphi_{\Phi (g,3)(\mjs_q,1)}=((\mj^{(j-1)
m},\mj^i,\tau',\mj^{(l-j-1)m-1},\tau',\vmj)^q)=\varphi_{\Phi'
(g,3)(\mjs_q,1)}$.
\end{tabbing}

If $x=1$, $y=3$, and $z=2$, we modify all the cases of Lemma 4.3
so that $\kappa$, $\beta$, $\tau$, and $\iota$ are replaced by
$\kappa''$, $\beta''$, $\tau''$, and $\iota''$, and the tuples
representing the natural transformations are multiplied by the
matrix $\displaystyle I_{n'}\otimes (\underbrace{1,\ldots,1}_q)
\otimes I_{m'}$. For example, Case 1.1.1 now reads
\begin{tabbing}
\hspace{1.5em}\=$\varphi_{\Phi (g,3)(\mjs_q,2)}$\=$=(\mj^{(j-1)
mq},(\mj^i,\tau'',\mj^{m-i-1})^q,\mj^{(l-j-2)mq},(\mj^i,\tau'',\mj^{m-i-1})^q,\vmj)$
\\[1ex]
\>\>$=\varphi_{\Phi' (g,3)(\mjs_q,2)}$.\`$\dashv$
\end{tabbing}

By relying on Lemma 6.3, we can prove a lemma analogous to Lemma
4.4 and this suffices for the proof of $(i)$ by induction on the
sum of lengths of ${\rm H}_1$ and $\Gamma_1$. We can prove $(ii)$
in a dual manner. Hence, we have:

\prop{Theorem 6.5}{The three-fold reduced bar construction $\WM$,
together with the natural transformations $\omega$, makes a lax
functor from $(\Delta^{\!op})^3$ to \emph{Cat}.}

\section{The $n$-fold monoidal categories}

The notion of $n$-fold monoidal category that we use in this paper
is defined in \cite[Section 7.6]{AM10} under the name
$n$-\textit{monoidal category}. Before we define the notion of
$(n+1)$-fold monoidal category, for $n\geq 3$, we first define
what the arrows between the $n$-fold monoidal categories are. For
this inductive definition we assume that an $n$-fold monoidal
category, for $n\geq 3$, is a category $\cal M$ equipped with $n$
monoidal structures ${\cal M}_1=\langle{\cal M},\otm_1,
I_1\rangle,\ldots,{\cal M}_n=\langle{\cal M},\otm_n, I_n\rangle$
such that for every $1\leq k<l<m\leq n$, the category $\cal M$
with the structures ${\cal M}_k$, ${\cal M}_l$ and ${\cal M}_m$ is
three-fold monoidal. Hence, for every $1\leq k<l\leq n$, the
category $\cal M$ with the structures ${\cal M}_k$ and ${\cal
M}_l$ is two-fold monoidal. We denote by $\kappa_{k,l}\!: I_k\str
I_l$, $\beta_{k,l}\!: I_k\str I_k\otm_l I_k$, $\tau_{k,l}\!:
I_l\otm_k I_l\str I_l$ and
\[
\iota_{k,l}\!:(A\otm_l B)\otm_k(C\otm_l D)\str(A\otm_k
C)\otm_l(B\otm_k D)
\]
the required arrows of $\cal M$.

\vspace{2ex}

\noindent{\sc Definition.}\quad An $n$\emph{-fold monoidal
functor} between two $n$-fold monoidal categories $\cal C$ and
$\cal D$ is a $(2n+1)$-tuple $\langle
F,\sigma^1,\delta^1,\ldots,\sigma^n,\delta^n\rangle$, where for
$k\in\{1,\ldots,n\}$,
\[
\sigma^k_{A,B}\!:FA\otimes^{\cal D}_k FB\str F(A\otimes^{\cal C}_k
B)\quad{\rm and}\quad \delta^k\!:I^{\cal D}_k\str FI^{\cal C}_k
\]
are arrows of $\cal D$ natural in $A$ and $B$, such that $\langle
F,\sigma^k,\delta^k\rangle$ is a monoidal functor between the
$k$th monoidal structures of $\cal C$ and $\cal D$. Moreover, for
every $1\leq k<l\leq n$, the following four diagrams commute (with
the superscripts $\cal C$ and $\cal D$ omitted):

\begin{center}
\begin{picture}(80,50)
\put(0,40){\makebox(0,0){$I_k$}} \put(80,40){\makebox(0,0){$I_l$}}
\put(0,10){\makebox(0,0){$FI_k$}}
\put(80,10){\makebox(0,0){$FI_l$}}

\put(15,10){\vector(1,0){50}} \put(15,40){\vector(1,0){50}}
\put(0,33){\vector(0,-1){15}} \put(80,33){\vector(0,-1){15}}

\put(40,7){\makebox(0,0)[t]{$F\kappa_{k,l}$}}
\put(40,43){\makebox(0,0)[b]{$\kappa_{k,l}$}}
\put(-3,25){\makebox(0,0)[r]{$\delta^k$}}
\put(83,25){\makebox(0,0)[l]{$\delta^l$}}

\end{picture}
\end{center}

\begin{center}
\begin{picture}(280,80)
\put(0,70){\makebox(0,0){$I_k$}}
\put(80,70){\makebox(0,0){$I_k\otimes_l I_k$}}
\put(0,10){\makebox(0,0){$FI_k$}}
\put(80,40){\makebox(0,0){$FI_k\otimes_l FI_k$}}
\put(80,10){\makebox(0,0){$F(I_k\otimes_l I_k)$}}

\put(15,10){\vector(1,0){35}} \put(15,70){\vector(1,0){40}}
\put(0,63){\vector(0,-1){45}} \put(80,63){\vector(0,-1){15}}
\put(80,33){\vector(0,-1){15}}

\put(35,7){\makebox(0,0)[t]{$F\beta_{k,l}$}}
\put(35,73){\makebox(0,0)[b]{$\beta_{k,l}$}}
\put(-3,40){\makebox(0,0)[r]{$\delta^k$}}
\put(83,55){\makebox(0,0)[l]{$\delta^k\otimes_l\delta^k$}}
\put(83,25){\makebox(0,0)[l]{$\sigma^l$}}

\put(280,70){\makebox(0,0){$I_l$}}
\put(200,70){\makebox(0,0){$I_l\otimes_k I_l$}}
\put(280,10){\makebox(0,0){$FI_l$}}
\put(200,40){\makebox(0,0){$FI_l\otimes_k FI_l$}}
\put(200,10){\makebox(0,0){$F(I_l\otimes_k I_l)$}}

\put(225,10){\vector(1,0){35}} \put(225,70){\vector(1,0){40}}
\put(280,63){\vector(0,-1){45}} \put(200,63){\vector(0,-1){15}}
\put(200,33){\vector(0,-1){15}}

\put(245,7){\makebox(0,0)[t]{$F\tau_{k,l}$}}
\put(245,73){\makebox(0,0)[b]{$\tau_{k,l}$}}
\put(283,40){\makebox(0,0)[l]{$\delta^l$}}
\put(197,55){\makebox(0,0)[r]{$\delta^l\otimes_k\delta^l$}}
\put(197,25){\makebox(0,0)[r]{$\sigma^k$}}

\end{picture}
\end{center}

\begin{center}
\begin{picture}(200,80)

\put(0,70){\makebox(0,0){$(FA\otimes_l FB)\otimes_k(FC\otimes_l
FD)$}} \put(200,70){\makebox(0,0){$(FA\otimes_k
FC)\otimes_l(FB\otimes_k FD)$}}

\put(0,40){\makebox(0,0){$F(A\otimes_l B)\otimes_k F(C\otimes_l
D)$}} \put(200,40){\makebox(0,0){$F(A\otimes_k C)\otimes_l
F(B\otimes_k D)$}}

\put(0,10){\makebox(0,0){$F((A\otimes_l B)\otimes_k(C\otimes_l
D))$}} \put(200,10){\makebox(0,0){$F((A\otimes_k
C)\otimes_l(B\otimes_k D))$}}

\put(70,10){\vector(1,0){60}} \put(70,70){\vector(1,0){60}}
\put(0,63){\vector(0,-1){15}} \put(200,63){\vector(0,-1){15}}
\put(0,33){\vector(0,-1){15}} \put(200,33){\vector(0,-1){15}}

\put(100,7){\makebox(0,0)[t]{$F\iota_{k,l}$}}
\put(100,73){\makebox(0,0)[b]{$\iota_{k,l}$}}
\put(-3,55){\makebox(0,0)[r]{$\sigma^l\otimes_k\sigma^l$}}
\put(-3,25){\makebox(0,0)[r]{$\sigma^k$}}
\put(203,55){\makebox(0,0)[l]{$\sigma^k\otimes_l\sigma^k$}}
\put(203,25){\makebox(0,0)[l]{$\sigma^l$}}

\end{picture}
\end{center}

Let $Mon_n(Cat)$ be the 2-category whose 0-cells are the $n$-fold
monoidal categories, 1-cells are the $n$-fold monoidal functors,
and 2-cells are the $n$-\emph{fold monoidal transformations},
i.e.,\ monoidal transformations with respect to all $n$
structures. The monoidal structure of $Mon_n(Cat)$ is again given
by 2-products.

\vspace{2ex}

\noindent{\sc Definition.} An $(n+1)$\emph{-fold monoidal}
category is a pseudomonoid in $Mon_n(Cat)$.

\vspace{2ex}

By this inductive definition, it is clear that an $n$-fold
monoidal category satisfies the assumptions given above, which we
may take as an unfolded form of this definition. As in the case of
two-fold and three-fold monoidal categories, we are only
interested in $n$\emph{-fold strict monoidal} categories, i.e.,\
when the structures ${\cal M}_1,\ldots,{\cal M}_n$ are strict
monoidal.


\vspace{2ex}

The $n$-fold monoidal categories defined in
\cite[Definition~1.7]{BFSV} are the $n$-fold strict monoidal
categories from above in which, moreover, it is assumed that
$I_1=\ldots=I_n=0$, and all the $\kappa$, $\beta$ and $\tau$
arrows are replaced by the identity $\mj_0$. Also, for every $n$,
a symmetric monoidal category is $n$-fold monoidal with all $n$
monoidal structures being the same.

On the other hand, it is not true that every $n$-fold strict
monoidal category in our sense is an $n$-fold monoidal in the
sense of \cite{BFSV}. It is not only the case that the difference
would appear in arrows that involve the units. The arrows of the
form
\[
A\otimes_i B\str A\otimes_j B\quad {\rm and}\quad A\otimes_i B\str
B\otimes_j A,
\]
for $i<j$ (see \cite[Remark~1.4]{BFSV}), show that the
axiomatization of $n$-fold monoidal categories given in
\cite{BFSV} leads to a non-conservative extension of its fragment
without units. These arrows are not presumed by our definition.
Hence, the categories would be different in their unit-free
fragments too.


\section{The $n$-fold reduced bar construction}


In Sections~4 and~6, we have defined the $n$-fold reduced bar
construction for $n=2$ and $n=3$. We define, in the same manner,
the $n$-fold reduced bar construction for arbitrary $n\geq 3$.
This construction corresponds to the one given in the proof of
\cite[Theorem~2.1]{BFSV}, save that the latter construction is
based on a category that is $n$-fold monoidal in the sense of that
paper.

For an $n$-fold strict monoidal category $\cal M$, we define
functions $\WM$ from objects and arrows of $(\Delta^{\!op})^n$ to
objects and arrows of \emph{Cat} in the following manner.

\vspace{2ex}

\noindent {\sc Definition.}\quad The \emph{$n$-fold reduced bar
construction} $\WM$ is defined on objects of $(\Delta^{\!op})^n$
as:
\[
\WM(k_1,\ldots,k_n)=\cM^{k_1\cdot \ldots\cdot k_n},
\]
and for arrows $f_k:s_k\str t_k$, $1\leq k\leq n$, of
$\Delta^{\!op}$, we define $\WM(f_1,\ldots,f_n)$ as the
composition
\[
(\WM_n(f_n))^{t_1\cdot\ldots\cdot
t_{n-1}}\circ\ldots\circ(\WM_k^{s_{k+1}\cdot\ldots\cdot
s_n}(f_k))^{t_1\cdot\ldots\cdot t_{k-1}}\circ\ldots\circ
\WM_1^{s_2\cdot\ldots\cdot s_n}(f_1).
\]

For example, for $\cal M$ being a four-fold strict monoidal
category, the functor $\WM(d^2_1,s^2_1,d^2_1,s^2_0):{\cal M}^4\str
{\cal M}^4$ is defined as the composition
\[
(\WM_4(s^2_0))^2(\WM_3(d^2_1))^2\circ \WM_2^{2}(s^2_1)\circ
\WM_1^2(d^2_1),
\]
and for an object $(A,B,C,D)$ of ${\cal M}^4$ we have that
\[
\WM(d^2_1,s^2_1,d^2_1,s^2_0)(A,B,C,D)=(I_4,(A\otimes_1
C)\otimes_3(B\otimes_1 D),I_4,I_2\otimes_3 I_2).
\]

In order to prove that $\WM$ is a lax functor, for every
composable pair of arrows $e_1$ and $e_2$ of $(\Delta^{\!op})^n$,
we have to define a natural transformation
\[
\omega_{e_2,e_1}\!:\WM(e_2)\circ \WM(e_1)\strt \WM(e_2\circ e_1),
\]
such that Diag~4.1 commutes. For this we use again coloured
sequences and their shuffles.

Let $\Phi_1,\ldots,\Phi_n$ be sequences of colours $1,\ldots,n$,
respectively and let $\Theta$ be a shuffle of these $n$ sequences.
For every member $(f,k)$ of $\Theta$, we define its \emph{inner
power} and its \emph{outer power} to be
\[
\prod_{k<l\leq n}t_l\quad {\rm and}\quad \prod_{1\leq l<k}t_l,
\]
respectively, where $t_l$ is the target of its right-closest
member of $\Theta$ of colour $l$ (again with adding appropriate
identities if necessary). We assume that the empty product is~1.
This definition is in accordance with the corresponding
definitions for two and three-fold cases; the difference is that
the powers fixed to be~1 (like, for example, the outer power of
$(f,1)$) are not mentioned there.

Let $\cal M$ be an $n$-fold strict monoidal category and let our
sequences be such that for every $1\leq k\leq n$,
$\circ\Phi_k\!:s_k\str t_k$. We define a functor
\[
\WM_\Theta\!:{\cal M}^{s_1\cdot\ldots\cdot s_n}\str{\cal
M}^{t_1\cdot\ldots\cdot t_n}
\]
in the following way: replace in $\Theta$ every $(f,k)$ whose
inner power is $i$ and outer power is $o$ by $(\WM_k^i(f))^o$, and
insert $\circ$'s.

It is easy to see that for arrows $f_k$, $1\leq k\leq n$, of
$\Delta^{\!op}$, we have that
\[
\WM(f_1,\ldots,f_n)=\WM_{\Phi_n\ldots\Phi_1},
\]
for arbitrary sequences $\Phi_k$ of colour~$k$, $1\leq k\leq n$,
such that $\circ\Phi_k=f_k$. This may serve as an alternative
(\emph{combinatorial}) definition of the $n$-fold reduced bar
construction $\WM$.

We define the natural transformations $\omega$ following the lines
of Sections 4 and 6. In order to compare some notions needed for
this definition with the corresponding notions introduced in
Sections~4 and~6, we use the symbol $n$ for an object of
$\Delta^{\!op}$. To prevent ambiguities, we introduce a new symbol
$\dot{n}$, and assume that our category $\cal M$ is $\dot{n}$-fold
strict monoidal and that $\WM$ is the $\dot{n}$-fold reduced bar
construction. This includes just a few occurrences of $\dot{n}$
ending with Lemma~8.1, when we return to the standard notation.

For basic arrows $f\!:n\str n'$ and $g\!:m\str m'$ of
$\Delta^{\!op}$, for $k,l$ such that $0\leq k<l\leq \dot{n}$, and
$u,v,w\geq 0$, we define a natural transformation
\[
\chi^{k,l}_{u,v,w}(f,g)\!:(\WM_k^{vm'\!w}(f))^u\circ(\WM_l^w(g))^{unv}\;\strt\;
(\WM_l^w(g))^{un'\!v}\circ(\WM_k^{vmw}(f))^u
\]
to be the identity natural transformation except in the following
cases:

\begin{center}
\begin{tabular}{c|c|c}
$f$ & $g$ & $\chi^{k,l}_{u,v,w}(f,g)$
\\
\hline\hline\\[-2.2ex] $s^{n+1}_j$ & $s^{m+1}_i$ &
$\displaystyle((\mj^{j(m+1)vw},(\underbrace{\mj^{iw},\kappa^w_{k,l},\vmj}_{(m+1)w})^v,\vmj)^u)$
\\
\hline\\[-2.2ex] $d^n_j$, $1\leq j\leq n-1$ & $s^{m+1}_i$ &
$\displaystyle((\mj^{(j-1)(m+1)vw},(\underbrace{\mj^{iw},\tau^w_{k,l},\vmj}_{(m+1)w})^v,\vmj)^u)$
\\
\hline\\[-2.2ex] $s^{n+1}_j$ & $d^m_i$, $1\leq i\leq m-1$ &
$\displaystyle((\mj^{j(m-1)vw},(\underbrace{\mj^{(i-1)w},\beta^w_{k,l},\vmj}_{(m-1)w})^v,\vmj)^u)$
\\
\hline\\[-2.2ex] $d^n_j$, $1\leq j\leq n-1$ & $d^m_i$, $1\leq i\leq m-1$ &
$\displaystyle((\mj^{(j-1)(m-1)vw},(\underbrace{\mj^{(i-1)w},\iota^w_{k,l},\vmj}_{(m-1)w})^v,\vmj)^u)$

\end{tabular}
\end{center}

Note that the tuples of the third column of the table above are
obtained as a result of multiplication of the tuples in the third
column of Table~\ref{table:2dim} (where $\kappa$, $\beta$, $\tau$,
and $\iota$ are replaced by $\kappa_{k,l}$, $\beta_{k,l}$,
$\tau_{k,l}$, and $\iota_{k,l}$) by the matrix
\[
(\underbrace{1,\ldots,1}_u)\otimes I_{n'} \otimes
(\underbrace{1,\ldots,1}_v) \otimes
I_{m'}\otimes(\underbrace{1,\ldots,1}_w).
\]

For the following lemma, which is analogous to Lemma 6.1, we
assume that $f\!:n\str n'$, $g\!:m\str m'$, and $h\!:p\str p'$ are
basic arrows of $\Delta^{\!op}$, that $1\leq a< b< c\leq \dot{n}$,
that $\Lambda$ is a shuffle of sequences of colours
$1,\ldots,\dot{n}$ with only identity arrows in it, and that
\[
u=\prod_{1\leq l< a}t_l,\quad v_1=\prod_{a<l<b}t_l,\quad
v_2=\prod_{b<l<c}t_l,\quad w=\prod_{b<l\leq \dot{n}}t_l,
\]
where $t_l$ is the target of the leftmost member of $\Lambda$ of
colour $l$.

For example, if $\dot{n}=7$, $a=2$, $b=4$, $c=5$, and
\[
\Lambda=(\mj_2,1)(\mj_n,2)(\mj_3,3)(\mj_m,4)
(\mj_p,5)(\mj_5,6)(\mj_4,7),
\]
then $u=2$, $v_1=3$, $v_2=1$, and $w=20$.

\prop{Lemma 8.1}{The following diagram commutes:}

\begin{center}
\begin{picture}(200,135)

\put(100,130){\makebox(0,0){$\WM_{(f,a)(g,b)(h,c)\Lambda}$}}
\put(30,90){\makebox(0,0){$\WM_{(g,b)(f,a)(h,c)\Lambda}$}}
\put(170,90){\makebox(0,0){$\WM_{(f,a)(h,c)(g,b)\Lambda}$}}
\put(30,50){\makebox(0,0){$\WM_{(g,b)(h,c)(f,a)\Lambda}$}}
\put(170,50){\makebox(0,0){$\WM_{(h,c)(f,a)(g,b)\Lambda}$}}
\put(100,10){\makebox(0,0){$\WM_{(h,c)(g,b)(f,a)\Lambda}$}}

\put(70,120){\vector(-3,-2){30}} \put(130,120){\vector(3,-2){30}}
\put(40,80){\vector(0,-1){23}} \put(160,80){\vector(0,-1){23}}
\put(40,40){\vector(3,-2){30}} \put(160,40){\vector(-3,-2){30}}

\put(45,115){\makebox(0,0)[r]{$\chi^{a,b}_{u,v_1,v_2
p'w}(f,g)\WM_{(h,c)\Lambda}$}}

\put(35,70){\makebox(0,0)[r]{
$\WM_{(g,b)\Lambda}\;\chi^{a,c}_{u,v_1 m v_2,w}(f,h)$}}

\put(45,25){\makebox(0,0)[r]{$\chi^{b,c}_{un'v_1,v_2,w}(g,h)\WM_{(f,a)\Lambda}$}}

\put(155,115){\makebox(0,0)[l]{
$\WM_{(f,a)\Lambda}\;\chi^{b,c}_{unv_1,v_2,w}(g,h)$}}

\put(165,70){\makebox(0,0)[l]{$\chi^{a,c}_{u,v_1m'v_2,w}(f,h)\WM_{(g,b)\Lambda}$}}

\put(155,25){\makebox(0,0)[l]{
$\WM_{(h,c)\Lambda}\;\chi^{a,b}_{u,v_1,v_2pw}(f,g)$}}

\end{picture}
\end{center}

\dkz The tuples representing the natural transformations of the
left-hand side and the right-hand side of this diagram are
obtained by multiplying the corresponding tuples of the diagram in
Lemma 6.1 (where $\kappa$, $\kappa'$, and $\kappa''$ are replaced
by $\kappa_{a,b}$, $\kappa_{b,c}$, and $\kappa_{a,c}$, etc.) by
the matrix
\[
(\underbrace{1,\ldots,1}_u)\otimes I_{n'} \otimes
(\underbrace{1,\ldots,1}_{v_1}) \otimes I_{m'}
(\underbrace{1,\ldots,1}_{v_2}) \otimes I_{p'}
\otimes(\underbrace{1,\ldots,1}_w).
\]
Hence, Lemma 6.1 directly implies this lemma. \qed

Let $\Theta_0,\ldots,\Theta_j$, for $j\geq 0$, be shuffles of
$\Phi_1,\ldots,\Phi_n$ such that $\Theta_0=\Theta$ and
$\Theta_j=\Phi_n\ldots\Phi_1$, and if $j>0$, then for every $0\leq
i\leq j-1$ we have that for some $1\leq k<l\leq n$,
$\Theta_i=\Pi(f,k)(g,l)\Lambda$ and
$\Theta_{i+1}=\Pi(g,l)(f,k)\Lambda$. We call
$\Theta_0,\ldots,\Theta_j$ a \emph{normalizing path} starting with
$\Theta$. Its \emph{length} is $j$ and Proposition 4.1 still
holds.

For $u$, $v$, and $w$ being respectively
\[
\prod_{1\leq z< k}t_z,\quad \prod_{k<z<l}t_z,\quad \prod_{l<z\leq
n}t_z,
\]
where $t_z$ is the target of the leftmost member of $\Lambda$ of
colour $z$, we have that
\[
\varphi_i=\WM_\Pi\;\chi^{k,l}_{u,v,w}(f,g)\;\WM_\Lambda,
\]
is a natural transformation from $\WM_{\Theta_i}$ to
$\WM_{\Theta_{i+1}}$. We define
$\varphi(\Theta_0,\ldots,\Theta_j)$ as in the two-fold case and
for $\Theta'_0,\ldots,\Theta'_j$ being another normalizing path
starting with $\Theta$, the following theorem is proved in the
same manner as Theorem 6.2, relying on Lemma 8.1 instead of Lemma
6.1.

\prop{Theorem 8.2}{
$\varphi(\Theta_0,\ldots,\Theta_j)=\varphi(\Theta'_0,\ldots,\Theta'_j)$.}

By Theorem 8.2, the following definition is correct.

\vspace{2ex}

\noindent {\sc Definition.}\quad Let $\varphi_\Theta:\WM_\Theta
\strt \WM_{\Phi_n\ldots\Phi_1}$ be
$\varphi(\Theta_0,\ldots,\Theta_j)$, for an arbitrary normalizing
path $\Theta_0,\ldots,\Theta_j$ starting with $\Theta$.

\vspace{2ex}

We are ready to define a natural transformation
\[
\omega_{e_2,e_1}\!:\WM(e_2)\circ \WM(e_1)\strt \WM(e_2\circ e_1),
\]
for every composable pair of arrows $e_1=(f^1_1,\ldots,f^1_n)$ and
$e_2=(f^2_1,\ldots,f^2_n)$ of $(\Delta^{\!op})^n$.

\vspace{2ex}

\noindent {\sc Definition.}\quad Let $\Phi^1_k$ and $\Phi^2_k$,
for $1\leq k\leq n$, be sequences of colour $k$, such that
$\circ\Phi^1_k$ is $(f^1_k)^{\rm nf}$ and $\circ\Phi^2_k$ is
$(f^2_k)^{\rm nf}$. We define
\[
\omega_{e_2,e_1}\quad {\rm as}\quad
\varphi_{\Phi^2_n\ldots\Phi^2_1\Phi^1_n\ldots\Phi^1_1}.
\]

It remains to prove that our Diag~4.1 commutes. Let
$e_1=(f^1_1,\ldots,f^1_n)$, $e_2=(f^2_1,\ldots,f^2_n)$ and
$e_3=(f^3_1,\ldots,f^3_n)$ be such that the composition $e_3\circ
e_2\circ e_1$ is defined in $(\Delta^{\!op})^n$. Let $\Phi^1_k$,
$\Phi^2_k$ and $\Phi^3_k$, for $1\leq k\leq n$, be sequences of
colour $k$, such that $\circ\Phi^1_k$ is $(f^1_k)^{\rm nf}$,
$\circ\Phi^2_k$ is $(f^2_k)^{\rm nf}$ and $\circ\Phi^3_k$ is
$(f^3_k)^{\rm nf}$. In this case, Diag~4.1 reads

\begin{center}
\begin{picture}(100,75)(0,5)

\put(50,75){\makebox(0,0){$\WM(e_3)\circ\WM(e_2)\circ \WM(e_1)$}}
\put(-40,40){\makebox(0,0){$\WM(e_3\circ e_2)\circ \WM(e_1)$}}
\put(140,40){\makebox(0,0){$\WM(e_3)\circ\WM(e_2\circ e_1)$}}
\put(50,10){\makebox(0,0){$\WM(e_3\circ e_2\circ e_1)$}}

\put(-10,60){\vector(-2,-1){20}} \put(110,60){\vector(2,-1){20}}
\put(-30,30){\vector(2,-1){30}} \put(130,30){\vector(-2,-1){30}}

\put(-23,60){\makebox(0,0)[r]{\scriptsize
$\varphi_{\Phi^3_n\ldots\Phi^3_1\Phi^2_n\ldots\Phi^2_1}
\WM_{\Phi^1_n\ldots\Phi^1_1}$}}
\put(117,60){\makebox(0,0)[l]{\scriptsize
$\WM_{\Phi^3_n\ldots\Phi^3_1}\,
\varphi_{\Phi^2_n\ldots\Phi^2_1\Phi^1_n\ldots\Phi^1_1}$}}

\put(-20,20){\makebox(0,0)[r]{\scriptsize
$\varphi_{\Phi_n'\ldots\Phi_1'\Phi^1_n\ldots\Phi^1_1}$}}
\put(120,20){\makebox(0,0)[l]{\scriptsize
$\varphi_{\Phi^3_n\ldots\Phi^3_1\Phi_n''\ldots\Phi_1''}$}}

\end{picture}
\end{center}
where $\circ\Phi_k'$ is $(f^3_k\circ f^2_k)^{\rm nf}$ and
$\circ\Phi_k''$ is $(f^2_k\circ f^1_k)^{\rm nf}$.

\vspace{2ex}

By Theorem 8.2 we have the following commutative diagram

\begin{center}
\begin{picture}(100,75)(0,5)

\put(50,75){\makebox(0,0){$\WM(e_3)\circ\WM(e_2)\circ \WM(e_1)$}}
\put(-40,40){\makebox(0,0){$\WM(e_3\circ e_2)\circ \WM(e_1)$}}
\put(140,40){\makebox(0,0){$\WM(e_3)\circ\WM(e_2\circ e_1)$}}
\put(50,10){\makebox(0,0){$\WM(e_3\circ e_2\circ e_1)$}}

\put(-10,60){\vector(-2,-1){20}} \put(110,60){\vector(2,-1){20}}
\put(-30,30){\vector(2,-1){30}} \put(130,30){\vector(-2,-1){30}}
\put(50,65){\line(0,-1){5}} \put(50,50){\vector(0,-1){30}}

\put(-23,60){\makebox(0,0)[r]{\scriptsize
$\varphi_{\Phi^3_n\ldots\Phi^3_1\Phi^2_n\ldots\Phi^2_1}
\WM_{\Phi^1_n\ldots\Phi^1_1}$}}
\put(117,60){\makebox(0,0)[l]{\scriptsize
$\WM_{\Phi^3_n\ldots\Phi^3_1}\,
\varphi_{\Phi^2_n\ldots\Phi^2_1\Phi^1_n\ldots\Phi^1_1}$}}

\put(-20,20){\makebox(0,0)[r]{\scriptsize
$\varphi_{\Phi^3_n\Phi^2_n\ldots\Phi^3_1\Phi^2_1\Phi^1_n
\ldots\Phi^1_1}$}} \put(120,20){\makebox(0,0)[l]{\scriptsize
$\varphi_{\Phi^3_n\ldots\Phi^3_1\Phi^2_n\Phi^1_n\ldots
\Phi^2_1\Phi^1_1}$}} \put(50,55){\makebox(0,0){\scriptsize
$\varphi_{\Phi^3_n\ldots\Phi^3_1\Phi^2_n\ldots\Phi^2_1\Phi^1_n
\ldots\Phi^1_1}$}}

\end{picture}
\end{center}
Hence, to prove that Diag~4.1 commutes, it suffices to show that

\begin{tabbing}
\hspace{1.5em}\=$(i)$\quad\=$\varphi_{\Phi^3_n\Phi^2_n\ldots\Phi^3_1\Phi^2_1\Phi^1_n
\ldots\Phi^1_1}=
\varphi_{\Phi_n'\ldots\Phi_1'\Phi^1_n\ldots\Phi^1_1}$\quad and
\\[2ex]
\>$(ii)$\>$\varphi_{\Phi^3_n\ldots\Phi^3_1\Phi^2_n\Phi^1_n\ldots
\Phi^2_1\Phi^1_1}=
\varphi_{\Phi^3_n\ldots\Phi^3_1\Phi_n''\ldots\Phi_1''}$.
\end{tabbing}

To prove $(i)$ and $(ii)$ we use the same arguments as in the
two-fold case. Let $\Lambda$ be a shuffle of sequences of colours
$1,\ldots,n$ with only identity arrows in it. Let $1\leq k<l\leq
n$ and let $u$, $v$, and $w$ be respectively
\[
\prod_{1\leq z< k}t_z,\quad \prod_{k<z<l}t_z,\quad \prod_{l<z\leq
n}t_z,
\]
where $t_z$ is the target of the leftmost member of $\Lambda$ of
colour $z$.

\prop{Lemma 8.3}{If $\Phi$ and $\Phi'$ are sequences of colour $k$
such that $\circ\Phi=\circ\Phi'$ is a basic equation of
$\Delta^{\!op}$, and $g$ is a basic arrow of $\Delta^{\!op}$, then
we have that $\varphi_{\Phi (g,l) \Lambda}=\varphi_{\Phi'
(g,l)\Lambda}$.}

\dkz Suppose the target of $\circ\Phi$ is $n'$ and the target of
$g$ is $m'$. We proceed as in Lemma 4.3 with all the cases
modified so that $\kappa$, $\beta$, $\tau$, and $\iota$ are
replaced by $\kappa_{k,l}$, $\beta_{k,l}$, $\tau_{k,l}$, and
$\iota_{k,l}$, and the tuples representing the natural
transformations are multiplied by the matrix
\[
(\underbrace{1,\ldots,1}_u)\otimes I_{n'}\otimes
(\underbrace{1,\ldots,1}_v)\otimes
I_{m'}\otimes(\underbrace{1,\ldots,1}_w).
\]
So, for example, Case 1.1.1 now reads
\begin{tabbing}
\hspace{1.5em}\=$\varphi_{\Phi (g,l)\Lambda}$\=$=$
\\[1ex]
\`$((\mj^{(j-1)
mvw},(\mj^{iw},\tau_{k,l}^w,\mj^{(m-i-1)w})^v,\mj^{(l-j-2)mvw},
(\mj^{iw},\tau_{k,l}^w,\mj^{(m-i-1)w})^v,\vmj)^u)$
\\[1ex]
\>\>$=\varphi_{\Phi' (g,l)\Lambda}$.\`$\dashv$
\end{tabbing}

By relying on Lemma 8.3, we can prove a lemma analogous to Lemma
4.4 and this suffices for the proof of $(i)$ by induction on the
sum of lengths of $\Phi^1_n,\ldots,\Phi^1_2$. We can prove $(ii)$
in a dual manner. So, for every $n\geq 2$ we have:

\prop{Theorem 8.5}{The $n$-fold reduced bar construction $\WM$,
together with the natural transformations $\omega$, makes a lax
functor from $(\Delta^{\!op})^n$ to \emph{Cat}.}

We see, by analyzing this result, that the conditions imposed by
the definition of $n$-fold monoidal categories are not only
sufficient, but they are also necessary to prove the correctness
of the $n$-fold reduced bar construction. If one proves this
through the steps established by our Theorem 8.2 and Lemmata
analogous to Lemma 8.3, then all the combinatorial structure of
$n$-fold monoidal categories is used.

Since every $n$-fold monoidal category in the sense of \cite{BFSV}
is an $n$-fold strict monoidal category in our sense, Theorem~8.5
gives an alternative proof for \cite[Theorem~2.1]{BFSV}. Every
braided strict monoidal category is a two-fold monoidal category
in the sense of \cite{BFSV} and every symmetric strict monoidal
category is an $\infty$-monoidal category in the sense of
\cite{BFSV}. Hence, our Theorem~8.5 covers all the related results
concerning these categories. Also, the correctness of the reduced
bar construction of \cite[Lemma~7.1]{PT13} follows from this
theorem.

\section{Delooping}

This section, which is inspired by \cite[Section~2]{BFSV},
explains how to use Theorem~8.5 for delooping of classifying
spaces of $n$-fold monoidal categories. Theorem~2.2 of \cite{BFSV}
says that the group completion of the nerve of an $n$-fold
monoidal category is an $n$-fold loop space. It is an easy
corollary of a generalization of \cite[Proposition~1.5]{S74} and
\cite[Theorem~2.1]{BFSV}.

A formulation of a generalization of \cite[Proposition~1.5]{S74}
is given in \cite[paragraph preceding Theorem~2.1]{BFSV}. This
seems to be a folklore result amongst the experts, but we couldn't
find written proof, or a precise formulation of it. The note
\cite{P14} is prepared to rectify that. We sketch a delooping
procedure based on the results of this note.

For $m\geq 1$, consider the arrows $i_1,\ldots, i_m:m\str 1$ of
$\Delta^{op}$ given by the following diagrams.

\begin{center}
\begin{picture}(320,40)(30,0)

\put(40,20){\makebox(0,0)[r]{$i_1\!:$}}

\put(50,10){\circle*{2}} \put(70,10){\circle*{2}}
\put(50,30){\circle*{2}} \put(70,30){\circle*{2}}
\put(85,30){\makebox(0,0){\ldots}} \put(100,30){\circle*{2}}

\put(50,2){\makebox(0,0)[b]{\scriptsize $0$}}
\put(70,2){\makebox(0,0)[b]{\scriptsize $1$}}
\put(50,33){\makebox(0,0)[b]{\scriptsize $0$}}
\put(70,33){\makebox(0,0)[b]{\scriptsize $1$}}
\put(100,33){\makebox(0,0)[b]{\scriptsize $m$}}

\put(50,10){\line(0,1){20}} \put(70,10){\line(0,1){20}}

\put(140,20){\makebox(0,0)[r]{$i_2\!:$}}

\put(150,10){\circle*{2}} \put(170,10){\circle*{2}}
\put(150,30){\circle*{2}} \put(170,30){\circle*{2}}
\put(190,30){\circle*{2}} \put(205,30){\makebox(0,0){\ldots}}
\put(220,30){\circle*{2}}

\put(150,2){\makebox(0,0)[b]{\scriptsize $0$}}
\put(170,2){\makebox(0,0)[b]{\scriptsize $1$}}
\put(150,33){\makebox(0,0)[b]{\scriptsize $0$}}
\put(170,33){\makebox(0,0)[b]{\scriptsize $1$}}
\put(190,33){\makebox(0,0)[b]{\scriptsize $2$}}
\put(220,33){\makebox(0,0)[b]{\scriptsize $m$}}

\put(150,10){\line(1,1){20}} \put(170,10){\line(1,1){20}}

\put(250,20){\makebox(0,0){\ldots}}

\put(290,20){\makebox(0,0)[r]{$i_m\!:$}}

\put(300,10){\circle*{2}} \put(320,10){\circle*{2}}
\put(300,30){\circle*{2}} \put(330,30){\circle*{2}}
\put(315,30){\makebox(0,0){\ldots}} \put(350,30){\circle*{2}}

\put(300,2){\makebox(0,0)[b]{\scriptsize $0$}}
\put(320,2){\makebox(0,0)[b]{\scriptsize $1$}}
\put(300,33){\makebox(0,0)[b]{\scriptsize $0$}}
\put(330,33){\makebox(0,0)[b]{\scriptsize $m\!-\!1$}}
\put(350,33){\makebox(0,0)[b]{\scriptsize $m$}}

\put(300,10){\line(3,2){30}} \put(320,10){\line(3,2){30}}

\end{picture}
\end{center}
These arrows are related to projections, which is explained in
\cite[Section~2]{P15} and \cite[Section~3]{P14}.

We use the following notation in the sequel. For functors
$F_i\!:{\cal A}\str {\cal B}_i$, $1\leq i\leq m$, let $\langle
F_1,\ldots,F_m\rangle: {\cal A}\str {\cal B}_1\times\ldots\times
{\cal B}_m$ be the functor obtained by the Cartesian structure of
\emph{Cat}.

Let $\WM$ be the $n$-fold reduced bar construction for $n\geq 2$.
It is easy to verify that for every $l\in\{0,\ldots, n-1\}$ and
every $k\geq 0$, the functor
$W\!:\Delta^{op}\str\mbox{\emph{Cat}}$ defined as
\[
\WM(\underbrace{1,\ldots,
1}_{l},\underline{\hspace{1ex}}\,,k,\ldots, k)
\]
is such that
\[
\langle W(i_1),\ldots,W(i_m) \rangle\!: W(m)\str(W(1))^m
\]
is the identity. This means that $\WM$ is Segal's lax functor
according to \cite[Definition~4.2]{P14}.

Let $V$ be a rectification of $\WM$ obtained by
\cite[Theorem~2]{S72}, and let \linebreak
${B\!:\mbox{\it{Cat}}\str\mbox{\it{Top}}}$ be the
\emph{classifying space} functor, i.e., the composition
$|\mbox{\hspace{1ex}}|\circ N$, where $N\!:\mbox{\it{Cat}}\str
\mbox{\it{Top}}^{\Delta^{op}}$ is the \emph{nerve} functor, and
$|\mbox{\hspace{1ex}}|:\mbox{\it{Top}}^{\Delta^{op}}\str
\mbox{\it{Top}}$ is the standard \emph{geometric realization}
functor. By \cite[Corollary~4.4]{P14}, $B\circ V$ is a
multisimplicial space such that for $X$ being the simplicial space
defined as
\[
(B\circ V)(\underbrace{1,\ldots,
1}_{l},\underline{\hspace{1ex}}\,,k,\ldots, k),
\]
the map
\[
\langle X(i_1),\ldots,X(i_m) \rangle\!: X(m)\str(X(1))^m
\]
is a homotopy equivalence.

By applying \cite[Lemma~3.1]{P14} to the simplicial space $(B\circ
V)(1,\ldots, 1,\underline{\hspace{1ex}}\,)$, we obtain a homotopy
associative H-space structure on $(B\circ V)(1,\ldots, 1)$. The
following theorem (in which $|\mbox{\hspace{1ex}}|$ denotes the
standard geometric realization of multisimplicial spaces) is taken
over from \cite[Theorem~5.1]{P14}.

\prop{Theorem 8.6}{If $(B\circ V)(1,\ldots, 1)$, with respect to
the above \emph{H}-space structure is grouplike, then $B{\cal
M}\simeq \Omega^n|B\circ V|$. }

Hence, up to group completion, the realization $|B\circ V|$ of the
multisimplicial space $B\circ V$ is an $n$-fold delooping of the
classifying space $B{\cal M}$ of $\cal M$.

\vspace{2ex}

\emph{Bicartesian categories}, i.e.,\ categories with all finite
coproducts and products may serve as examples of $n$-fold monoidal
categories that are not $n$-fold monoidal in the sense of
\cite{BFSV}. If we denote the nullary and binary coproducts of a
bicartesian category by $0$ and $+$, and nullary and binary
products by $1$ and $\times$, then the unique arrows
\[
\kappa\!:0\str 1,\quad \beta\!:0\str 0\times 0,\quad
\tau\!:1+1\str 1
\]
of this category together with the arrows
\[
\iota_{A,B,C,D}\!:(A\times B)+(C\times D)\str (A+C)\times(B+D),
\]
which are canonical in the coproduct-product structures, guarantee
that such a category may be conceived as a two-fold monoidal with
the first monoidal structure given by $+$ and $0$, and the second
given by $\times$ and $1$. Furthermore, such a category may be
conceived as an $n$-fold monoidal category in $n+1$ different ways
by taking first $0\leq k\leq n$ monoidal structures to be given by
the symmetric monoidal structure brought by $+$ and $0$, and the
remaining $n-k$ monoidal structures to be given by the symmetric
monoidal structure brought by $\times$ and $1$.

As a consequence of this fact there is a family, indexed by pairs
of natural numbers, of reduced bar constructions based on a
bicartesian category (strictified in both monoidal structures).
This is related to Adams' remark on $E_\infty$ ring spaces given
in \cite[\S 2.7]{A78} where the bicartesian category \emph{FinSet}
of finite sets and functions, with disjoint union as $+$ and
Cartesian product as $\times$, is mentioned. According to Segal,
\cite[\S 2]{S74},``most fundamental $\Gamma$-space'' arises from
this category under disjoint union.

By applying our results, it is possible to combine the disjoint
union and Cartesian product in the category \emph{FinSet} to
obtain various multisimplicial spaces. Since we have the initial
(and a terminal) object in \emph{FinSet}, its classifying space is
contractible and all the other realizations of simplicial sets in
question are path-connected. Hence, the induced H-space structures
are grouplike, and there is no need for group completion when one
starts to deloop \emph{FinSet} with respect to the disjoint union
and then continue to deloop it with respect to Cartesian product.
However, all these deloopings are contractible.

Since the notion of $n$-fold monoidal category is equationally
presented, there are $n$-fold monoidal categories freely generated
by sets of objects. We believe that delooping of classifying
spaces of such categories deserves particular attention. Also,
some other examples of $n$-fold monoidal categories from the
literature (e.g. \cite[Sections~6.4 and 7.3]{AM10}) could be
interesting from the point of view of delooping.


\section{Appendix}

By the definition given in Section~2, a two-fold monoidal category
consists of the following:

\vspace{.5ex}

1. a monoidal category $\langle{\cal M},\otm_1, I_1,
\alpha_1,\rho_1,\lambda_1\rangle$ (here $\alpha_1$, $\rho_1$, and
$\lambda_1$, respectively, denote associativity, right and left
identity natural isomorphisms),

\vspace{.5ex}

2. monoidal functors $\otm_2\!:{\cal M}\times {\cal M}\str {\cal
M}$ and $ I_2\!: 1\str {\cal M}$,

\vspace{.5ex}

3. monoidal transformations $\alpha_2$, $\rho_2$, and $\lambda_2$
such that $\langle{\cal M},\otm_2, I_2,
\alpha_2,\rho_2,\lambda_2\rangle$ satisfies the pseudomonoid
conditions (i.e.,\ the equations of a monoidal category).

\vspace{.3ex}

That $\otm_2$ is a monoidal functor means that there is a natural
transformation $\iota$ given by the family of arrows
\[
\iota_{A,B,C,D}\!:(A\otm_2 B)\otm_1(C\otm_2 D)\str(A\otm_1
C)\otm_2(B\otm_1 D),
\]
and an arrow $\beta\!: I_1\str  I_1\otm_2 I_1$ such that the
following three diagrams commute:

\vspace{3ex}

\begin{center}
\begin{picture}(200,75)
\small \put(0,10){\makebox(0,0){$(A\otm_2 D)\otm_1((B\otm_2
E)\otm_1(C\otm_2 F))$}} \put(200,10){\makebox(0,0){$((A\otm_2
D)\otm_1(B\otm_2 E))\otm_1(C\otm_2 F)$}}

\put(0,40){\makebox(0,0){$(A\otm_2 D)\otm_1((B\otm_1
C)\otm_2(E\otm_1 F))$}} \put(200,40){\makebox(0,0){$((A\otm_1
B)\otm_2(D\otm_1 E))\otm_1(C\otm_2 F)$}}

\put(0,70){\makebox(0,0){$(A\otm_1(B\otm_1
C))\otm_2(D\otm_1(E\otm_1 F))$}}
\put(200,70){\makebox(0,0){$((A\otm_1 B)\otm_1 C)\otm_2((D\otm_1
E)\otm_1 F)$}}

\put(85,10){\vector(1,0){27}} \put(85,70){\vector(1,0){27}}
\put(0,18){\vector(0,1){14}} \put(200,18){\vector(0,1){14}}
\put(0,48){\vector(0,1){14}} \put(200,48){\vector(0,1){14}}

\put(100,3){\makebox(0,0)[t]{$\alpha_1$}}
\put(100,77){\makebox(0,0)[b]{$\alpha_1\otm_2\alpha_1$}}
\put(-5,25){\makebox(0,0)[r]{$\mj\otm_1\iota$}}
\put(-5,55){\makebox(0,0)[r]{$\iota$}}
\put(205,25){\makebox(0,0)[l]{$\iota\otm_1\mj$}}
\put(205,55){\makebox(0,0)[l]{$\iota$}}

\put(100,40){\makebox(0,0){\small $(1)$}}

\end{picture}
\end{center}
\normalfont

\vspace{3ex}

\begin{center}
\begin{picture}(130,45)
\put(0,10){\makebox(0,0){$(A\otm_2 B)\otm_1 I_1$}}
\put(130,10){\makebox(0,0){$A\otm_2 B$}}

\put(0,40){\makebox(0,0){$(A\otm_2 B)\otm_1( I_1\otm_2 I_1)$}}
\put(130,40){\makebox(0,0){$(A\otm_1 I_1)\otm_2(B\otm_1 I_1)$}}

\put(35,10){\vector(1,0){70}} \put(55,40){\vector(1,0){22}}
\put(0,18){\vector(0,1){14}} \put(130,32){\vector(0,-1){14}}

\put(70,7){\makebox(0,0)[t]{$\rho_1$}}
\put(65,44){\makebox(0,0)[b]{$\iota$}}
\put(-5,25){\makebox(0,0)[r]{$\mj\otm_1\beta$}}
\put(135,25){\makebox(0,0)[l]{$\rho_1\otm_2\rho_1$}}

\put(65,25){\makebox(0,0){\small $(2)$}}

\end{picture}
\end{center}

\vspace{1ex}

\begin{center}
\begin{picture}(130,45)
\put(0,10){\makebox(0,0){$ I_1\otm_1(A\otm_2 B)$}}
\put(130,10){\makebox(0,0){$A\otm_2 B$}}

\put(0,40){\makebox(0,0){$( I_1\otm_2 I_1)\otm_1(A\otm_2 B)$}}
\put(130,40){\makebox(0,0){$( I_1\otm_1 A)\otm_2( I_1\otm_1 B)$}}

\put(35,10){\vector(1,0){70}} \put(55,40){\vector(1,0){22}}
\put(0,18){\vector(0,1){14}} \put(130,32){\vector(0,-1){14}}

\put(70,7){\makebox(0,0)[t]{$\lambda_1$}}
\put(65,44){\makebox(0,0)[b]{$\iota$}}
\put(-5,25){\makebox(0,0)[r]{$\beta\otm_1\mj$}}
\put(135,25){\makebox(0,0)[l]{$\lambda_1\otm_2\lambda_1$}}

\put(65,25){\makebox(0,0){\small $(3)$}}

\end{picture}
\end{center}

That $ I_2$ is a monoidal functor means that there are arrows
$\tau\!: I_2\otm_1 I_2\str  I_2$ and $\kappa\!: I_1\str  I_2$ such
that the following diagrams commute:

\vspace{1ex}

\begin{center}
\begin{picture}(90,75)
\put(45,10){\makebox(0,0){$ I_2$}}

\put(0,40){\makebox(0,0){$ I_2\otm_1 I_2$}}
\put(90,40){\makebox(0,0){$ I_2\otm_1 I_2$}}

\put(-10,70){\makebox(0,0){$ I_2\otm_1( I_2\otm_1 I_2)$}}
\put(100,70){\makebox(0,0){$( I_2\otm_1 I_2)\otm_1 I_2$}}

\put(0,65){\vector(0,-1){18}} \put(90,65){\vector(0,-1){18}}

\put(30,70){\vector(1,0){30}} \put(10,32){\vector(2,-1){30}}
\put(80,32){\vector(-2,-1){30}}

\put(45,73){\makebox(0,0)[b]{$\alpha_1$}}
\put(-3,55){\makebox(0,0)[r]{$\mj\otm_1\tau$}}
\put(93,55){\makebox(0,0)[l]{$\tau\otm_1\mj$}}

\put(23,23){\makebox(0,0)[tr]{$\tau$}}
\put(69,23){\makebox(0,0)[tl]{$\tau$}}

\put(45,48){\makebox(0,0){\small $(4)$}}
\end{picture}
\end{center}

\begin{center}
\begin{picture}(230,45)
\put(45,10){\makebox(0,0){$ I_2$}}

\put(0,40){\makebox(0,0){$ I_2\otm_1 I_1$}}
\put(90,40){\makebox(0,0){$ I_2\otm_1 I_2$}}

\put(20,40){\vector(1,0){50}} \put(10,32){\vector(2,-1){30}}
\put(80,32){\vector(-2,-1){30}}

\put(45,43){\makebox(0,0)[b]{$\mj\otm_1\kappa$}}
\put(23,23){\makebox(0,0)[tr]{$\rho_1$}}
\put(69,23){\makebox(0,0)[tl]{$\tau$}}

\put(45,28){\makebox(0,0){\small $(5)$}}

\put(185,10){\makebox(0,0){$ I_2$}}

\put(140,40){\makebox(0,0){$ I_1\otm_1 I_2$}}
\put(230,40){\makebox(0,0){$ I_2\otm_1 I_2$}}

\put(160,40){\vector(1,0){50}} \put(150,32){\vector(2,-1){30}}
\put(220,32){\vector(-2,-1){30}}

\put(185,43){\makebox(0,0)[b]{$\kappa\otm_1\mj$}}
\put(163,24){\makebox(0,0)[tr]{$\lambda_1$}}
\put(209,23){\makebox(0,0)[tl]{$\tau$}}

\put(185,28){\makebox(0,0){\small $(6)$}}

\end{picture}
\end{center}

Note that the unusual numbering of the following diagrams is due
to our wish to dualize the first six diagrams in some way, which
can clearly be seen from the list of twelve equations at the end
of Section~2. That $\alpha_2$ is a monoidal transformation means
that the following diagrams commute:

\vspace{1ex}

\begin{center}
\begin{picture}(200,70)(0,10)
\put(0,10){\makebox(0,0){$(A\otm_1 D)\otm_2((B\otm_1
E)\otm_2(C\otm_1 F))$}} \put(200,10){\makebox(0,0){$((A\otm_1
D)\otm_2(B\otm_1 E))\otm_2(C\otm_1 F)$}}

\put(0,40){\makebox(0,0){$(A\otm_1 D)\otm_2((B\otm_2
C)\otm_1(E\otm_2 F))$}} \put(200,40){\makebox(0,0){$((A\otm_2
B)\otm_1(D\otm_2 E))\otm_2(C\otm_1 F)$}}

\put(0,70){\makebox(0,0){$(A\otm_2(B\otm_2
C))\otm_1(D\otm_2(E\otm_2 F))$}}
\put(200,70){\makebox(0,0){$((A\otm_2 B)\otm_2 C)\otm_1((D\otm_2
E)\otm_2 F)$}}

\put(85,10){\vector(1,0){27}} \put(85,70){\vector(1,0){27}}
\put(0,32){\vector(0,-1){14}} \put(200,32){\vector(0,-1){14}}
\put(0,62){\vector(0,-1){14}} \put(200,62){\vector(0,-1){14}}

\put(100,3){\makebox(0,0)[t]{$\alpha_2$}}
\put(100,77){\makebox(0,0)[b]{$\alpha_2\otm_1\alpha_2$}}
\put(-5,25){\makebox(0,0)[r]{$\mj\otm_2\iota$}}
\put(-5,55){\makebox(0,0)[r]{$\iota$}}
\put(205,25){\makebox(0,0)[l]{$\iota\otm_2\mj$}}
\put(205,55){\makebox(0,0)[l]{$\iota$}}

\put(100,40){\makebox(0,0){\small $(7)$}}

\end{picture}
\end{center}

\vspace{5ex}

\begin{center}
\begin{picture}(90,70)

\put(45,10){\makebox(0,0){$ I_1$}}

\put(-10,70){\makebox(0,0){$ I_1\otm_2( I_1\otm_2 I_1)$}}
\put(100,70){\makebox(0,0){$( I_1\otm_2 I_1)\otm_2 I_1$}}

\put(0,40){\makebox(0,0){$ I_1\otm_2 I_1$}}
\put(90,40){\makebox(0,0){$ I_1\otm_2 I_1$}}

\put(0,47){\vector(0,1){18}} \put(90,47){\vector(0,1){18}}

\put(30,70){\vector(1,0){30}} \put(40,17){\vector(-2,1){30}}
\put(50,17){\vector(2,1){30}}

\put(45,73){\makebox(0,0)[b]{$\alpha_2$}}
\put(-3,55){\makebox(0,0)[r]{$\mj\otm_2\beta$}}
\put(93,55){\makebox(0,0)[l]{$\beta\otm_2\mj$}}

\put(23,23){\makebox(0,0)[tr]{$\beta$}}
\put(69,23){\makebox(0,0)[tl]{$\beta$}}

\put(45,48){\makebox(0,0){\small $(10)$}}

\end{picture}
\end{center}

That $\rho_2$ is a monoidal transformation means that the
following diagrams commute:

\vspace{1ex}

\begin{center}
\begin{picture}(130,40)
\put(0,10){\makebox(0,0){$A\otm_1 B$}}
\put(130,10){\makebox(0,0){$(A\otm_1 B)\otm_2 I_2$}}

\put(0,40){\makebox(0,0){$(A\otm_2 I_2)\otm_1(B\otm_2 I_2)$}}
\put(130,40){\makebox(0,0){$(A\otm_1 B)\otm_2( I_2\otm_1 I_2)$}}

\put(90,10){\vector(-1,0){68}} \put(55,40){\vector(1,0){22}}
\put(0,32){\vector(0,-1){14}} \put(130,32){\vector(0,-1){14}}

\put(65,7){\makebox(0,0)[t]{$\rho_2$}}
\put(65,44){\makebox(0,0)[b]{$\iota$}}
\put(-5,25){\makebox(0,0)[r]{$\rho_2\otm_1\rho_2$}}
\put(135,25){\makebox(0,0)[l]{$\mj\otm_2\tau$}}

\put(65,25){\makebox(0,0){\small $(8)$}}

\end{picture}
\end{center}

\vspace{3ex}

\begin{center}
\begin{picture}(90,40)
\put(45,10){\makebox(0,0){$ I_1$}}

\put(00,40){\makebox(0,0){$ I_1\otm_2 I_1$}}
\put(90,40){\makebox(0,0){$ I_1\otm_2 I_2$}}

\put(20,40){\vector(1,0){50}} \put(40,17){\vector(-2,1){30}}
\put(80,32){\vector(-2,-1){30}}

\put(45,43){\makebox(0,0)[b]{$\mj\otm_2\kappa$}}
\put(23,23){\makebox(0,0)[tr]{$\beta$}}
\put(69,23){\makebox(0,0)[tl]{$\rho_2$}}

\put(45,28){\makebox(0,0){\small $(11)$}}

\end{picture}
\end{center}

Finally, that $\lambda_2$ is a monoidal transformation means that
the following diagrams commute:

\vspace{1ex}

\begin{center}
\begin{picture}(130,40)
\put(0,10){\makebox(0,0){$A\otm_1 B$}}
\put(130,10){\makebox(0,0){$ I_2\otm_2(A\otm_1 B)$}}

\put(0,40){\makebox(0,0){$( I_2\otm_2 A)\otm_1( I_2\otm_2 B)$}}
\put(130,40){\makebox(0,0){$( I_2\otm_1 I_2)\otm_2(A\otm_1 B)$}}

\put(90,10){\vector(-1,0){68}} \put(55,40){\vector(1,0){22}}
\put(0,32){\vector(0,-1){14}} \put(130,32){\vector(0,-1){14}}

\put(70,7){\makebox(0,0)[t]{$\lambda_2$}}
\put(65,44){\makebox(0,0)[b]{$\iota$}}
\put(-5,25){\makebox(0,0)[r]{$\lambda_2\otm_1\lambda_2$}}
\put(135,25){\makebox(0,0)[l]{$\tau\otm_2\mj$}}

\put(65,25){\makebox(0,0){\small $(9)$}}

\end{picture}
\end{center}

\vspace{1ex}

\begin{center}
\begin{picture}(90,50)
\put(45,10){\makebox(0,0){$ I_1$}}

\put(00,40){\makebox(0,0){$ I_1\otm_2 I_1$}}
\put(90,40){\makebox(0,0){$ I_2\otm_2 I_1$}}

\put(20,40){\vector(1,0){50}} \put(40,17){\vector(-2,1){30}}
\put(80,32){\vector(-2,-1){30}}

\put(45,43){\makebox(0,0)[b]{$\kappa\otm_2\mj$}}
\put(23,23){\makebox(0,0)[tr]{$\beta$}}
\put(69,23){\makebox(0,0)[tl]{$\lambda_2$}}

\put(45,28){\makebox(0,0){\small $(12)$}}

\end{picture}
\end{center}

\end{document}